\documentstyle[12pt]{article}
\newcommand{\ga}{\alpha}
\newcommand{\gga}{\gamma}
\newcommand{\gG}{\Gamma}
\newcommand{\gb}{\beta}
\newcommand{\gd}{\delta}
\newcommand{\gk}{\kappa}

\newcommand{\gep}{\varepsilon}
\newcommand{\gvp}{\varphi}

\newcommand{\gl}{\lambda}

\newcommand{\gt}{\theta}

\newcommand{\gw}{\omega}
\newcommand{\go}{\omega}

\newcommand{\gz}{\zeta}

\newcommand{\bE}{\hbox{\hbox{$E$} \kern -9pt \raise -10pt \hbox{$\tilde{
}$}\kern 6pt}}
\newcommand{\bC}{\hbox{\hbox{$C$} \kern -9pt \raise -10pt \hbox{$\tilde{
}$}\kern 6pt}}
\newcommand{\bP}{\hbox{\hbox{$P$} \kern -7pt \raise -10pt \hbox{$\tilde{
}$}\kern 3pt}}
\newcommand{\bp}{\underline{p}}
\newcommand{\bG}{\underline{G}}
\newcommand{\bB}{\hbox{\hbox{$B$} \kern -9pt \raise -10pt \hbox{$\tilde{
}$}\kern 6pt}}
\newcommand{\bD}{\hbox{\hbox{$D$} \kern -9pt \raise -10pt \hbox{$\tilde{
}$}\kern 6pt}}

\newcommand{\bN}{\hbox{\bf N}}

\newcommand{\CA}{{\cal A}}

\newcommand{\CP}{{\cal P}}

\newcommand{\CT}{{\cal T}}
\newcommand{\CU}{{\cal U}}

\newcommand{\barP}{\bar{P}}
\newcommand{\barp}{\bar{p}}
\newcommand{\barN}{\bar{N}}

\newcommand{\ol}{\overline}

\newcommand{\eqdf}{\stackrel{def}{=}}

\newcommand{\Raro}{\Rightarrow}
\newcommand{\rarrow}{\rightarrow}

\newcommand{\raro}{\rightarrow}

\newcommand{\lan}{$\langle$}
\newcommand{\ran}{$\rangle$}
\newcommand{\itm}[1]{\item[(#1)]}
\newcommand{\itms}[1]{\item[[#1\kern -5pt]]}
\newcommand{\bd}{\begin{description}}
\newcommand{\ed}{\end{description}}

\newcommand{\ben}{\begin{enumerate}}
\newcommand{\een}{\end{enumerate}}

\newcommand{\nseqn}[2]{#1,\ldots,#2}

\newcommand{\seqn}[2]{\langle#1,\ldots ,#2\rangle}
\newcommand{\sseqn}[2]{\langle#1\,|#2\rangle}
\newcommand{\ssseqn}[2]{\langle#1\,|#2\rangle}

\newcommand{\setn}[2]{$\{\,$#1$\:|\:$#2$\}$}
\newcommand{\setm}[2]{\{\ #1\:|\:#2\}}

\newcommand{\fsetn}[2]{\{\,#1,\ldots ,#2\}}

\newcommand{\pair}[2]{\langle#1 ,#2\rangle}

\newcommand{\trpl}[3]{\langle#1 ,#2 ,#3 \rangle}

\newcommand{\nin}{\noindent}

\newcommand{\A}{\forall}
\newcommand{\rest}{\vbox{\hbox{$\:\kern -2pt\mathbin{\vert\kern-3.1pt\lower-1pt
   \hbox{$\mathsurround=0pt\mathchar"0012$}\kern-4pt}\:$}}}

\newcommand{\R}{\hbox{{\sf I}\kern-.1500em \hbox{\sf R}}}
\newcommand{\Q}
   {\hbox{${\rm Q} \kern -7.5pt \raise 2pt \hbox{\tiny$|$}\kern 7.5pt$}}
\newcommand{\dP}
   {\hbox{${\cal P}$}}
\newcommand{\dS}
   {\hbox{${\cal S}$}}
\newcommand{\dD}
   {\hbox{${\rm D} \kern -7.5pt  {\rm D}\kern 7.5pt$}}
\newcommand{\dq}
   {\hbox{${\cal Q}$}}
\newcommand{\dqq}
   {\hbox{${\cal Q}$}}
\newcommand{\C}
   {\hbox{${\rm C} \kern -7.5pt \raise 2pt \hbox{\tiny$|$}\kern 7.5pt$}}
\newcommand{\Z}{\hbox{\sf Z\kern-0.720em\hbox{ Z}}}

\newcommand{\empt}{\emptyset}

\newcommand{\rskip}{\vskip 9pt\noindent}
\newcommand{\ipar}[1]{{#1\quad}}

  \newcommand{\ovl}{\overline}

  \newcommand{\SA}{{\sf A}}
  \newcommand{\SB}{{\sf B}}
  \newcommand{\SC}{{\sf C}}
  \newcommand{\ST}{{\sf T}}
  \newcommand{\SS}{{\sf S}}
  \newcommand{\SP}{{\sf P}}
  \newcommand{\IN}{\varepsilon}
 
  \newcommand{\spa}{{\sf sp}}
  \newcommand{\st}{{\sf su}}
  \newcommand{\sua}{{\sf su}}
 
  \newcommand{\force}{{\raise 2pt\hbox{\small{$|$}}\kern -3pt \vdash \kern 3pt}}
  \newcommand{\Pbforce}{\force_{\dP_{\gb}}}
  \newcommand{\Pbnforce}{\force_{\dP_{\gb(n)}}}
  
  \newcommand{\Pbforces}[1]{\force_{\dP_{\gb_{#1}}}}

  \newcommand{\level}{{\it level}}
  \newcommand{\name}{{\it name}}
  \newcommand{\last}{{\it last}}
  \newcommand{\height}{{\it height}}
  \newcommand{\length}{{\it length}}
  \newcommand{\extends}{{\it extends}}
  \newcommand{\domain}{{\it domain}}
  \newcommand{\dom}{{\it dom}}

  \newcommand{\Aron}{Aronszajn}
  \newcommand{\proof}{\par\noindent{\bf Proof:}\quad}

   \newcommand{\bone}{\bf \mbox{{\small $1$}}}
\newcommand{\dline}[1]{\| #1 \|}

\newcommand{\bigtimes}{{\mbox{\large $\times$}}}
\newcommand{\fnn}[3]{#1:#2 \raro #3}

\def\newtheorems{\newtheorem{theorem}{Theorem}[section]

                 \newtheorem{lemma}[theorem]{Lemma}
                 \newtheorem{claim}[theorem]{Claim}

                 }

\newenvironment{defn}{\addtocounter{theorem}{1}\medskip\noindent{\bf
Definition \thetheorem\ \nolinebreak }\rm}{}
 
\newtheorems

\newenvironment{ntheorem}[1]{ \noindent{\bf #1\ }\em}{}

\begin{document}
\baselineskip 18pt
\setcounter{section}{0}
\title{A $\Delta^2_2$ Well-Order of the Reals And Incompactness of
$L(Q^{MM})$}
\author{Uri Abraham  \\ Department of mathematics and computer science\\
Ben-Gurion University, Be'er Sheva, ISRAEL;
\and
Saharon Shelah \thanks{Supported by a grant from the Israeli Academy of
Sciences. Pub. 403} \\ Institute of Mathematics\\ The Hebrew University,
Jerusalem, ISRAEL.}
\date{8-16-91}
\maketitle
\begin{abstract}
A  forcing poset of size $2^{2^{\aleph_1}}$ which adds no new reals
is described and shown to provide a $\Delta^2_2$ definable well-order of
the reals (in fact, any given relation of the reals may be so encoded in
some generic extension).
 The encoding of this well-order is obtained by playing with
products of \Aron\ trees: Some products are special while other are Suslin
trees. 

The paper also deals with the Magidor-Malitz logic: it is consistent that
this logic is highly non compact.
\end{abstract}
\nin{\bf\large Preface}
\label{sP}
\medskip

This paper deals with three issues: the question of definable
well-orders of the reals, the compactness of the  Magidor-Malitz
logic and the forcing techniques for specializing \Aron\ trees without
addition of new reals. 

In the hope of attracting a wider readership, we have tried to make this
paper as self contained as possible; in some cases we have reproved known
results, or given informal descriptions to remind the reader of what he or
she probably knows. We could not do so for the theorem that the iteration
of $\dD$-complete proper forcing adds no reals, and the reader may wish to
consult chapter V of
 Shelah [1982], or the new edition [1992].
 Anyhow, we rely on this theorem only at one point.

The question of the existence of a definable well-order of the set
of reals, $\R$, with all of its variants, is central in set theory.
As a starting point for the particular question which is studied here, we
take the theorem of Shelah and Woodin [1990] by which from the existence
of a large cardinal (a supercompact and even much less) it follows
that there is no well order of $\R$ in $L(\R)$.

 Assuming that there exists a cardinal which is
simultaneously measurable and Woodin, Woodin [in preparation] 
has shown that: {\it If CH holds,
then every $\Sigma^2_1$ set of reals is determined. Hence
there is no $\Sigma^2_1$ well order of the reals. }

A $\Sigma^2_i$ formula is a formula over the structure $\langle \bN,
\CP(\bN),\CP(\CP(\bN)), \in \ldots \rangle$
of type $\exists X_1\subseteq\CP (\bN) \A
X_2\subseteq \CP (\bN)\ldots \gvp(X_1\ldots)$, where there are $i$
alternations of quantifiers,
 and in $\gvp$ all quantifications are over $\bN$ and
$\CP(\bN)$. Equivalently, $\CP(\bN)$ can be replaced by $\langle
H=H(\gw_1),\in \rangle$ the collection of all hereditarily countable sets;
this seems to be useful in applications. so a $\Sigma^2_2$ formula has the
form $\exists X_1 \subseteq H \A X_2 \subseteq H \gvp(X_1\ldots)$, where
$\gvp$ is first order over $H$ and the $X_i$'s are predicates (subsets of
$H$).

A natural question asked by Woodin is whether his theorem cited above could not
be generalized to exclude   $\Sigma^2_2$ well-order of $\R$: Perhaps
 CH and some large cardinal
may imply that there is no $\Sigma^2_2$ well-order of $\R$. We give
a negative answer by providing a forcing poset (of small size,
$2^{2^{\aleph_0}}$) which adds no reals and 
gives generic extensions in which there exists
a $\Sigma^2_2$ well-order of $\R$.
Since supposedly any large cardinal retains its largeness after a
``small" forcing extension, no large cardinal contradicts a
$\Sigma^2_2$ well order of $\R$. 
 
  Specifically, we are going to prove the following Main Theorem.
 
  \begin{ntheorem}{Theorem A} Assume $\diamond_{\gw_{1}}$. Let $P(x)$ be a predicate
  (symbol). There is a (finite) sentence $\psi$ in the language containing
  $P(x)$ with the Magidor-Malitz quantifiers, such that the following holds.
  Given any $\SP\subseteq \gw_1$, (1) there is a model $M$ of $\psi$, enriching
  $(\gw_1,<,\SP)$ such that $P^M=\SP$, and (2) assuming
$2^{\aleph_0}=\aleph_1$ and  $2^{\aleph_1}=\aleph_2$ 
there is a forcing poset $Q$ of
  size $\aleph_2$ satisfying the $\aleph_2$-c.c. and adding no reals such
  that in $V^Q$
  \begin{tabular}{c}
  $M$ is the single model of $\psi$ (up to isomorphism).
  \end{tabular}
\end{ntheorem}
 
  Recall that the Magidor-Malitz logic $L(Q^{MM})$
 is obtained by adjoining to the regular
  first-order logic the quantifiers $Qxy\gvp(x,y)$ which is true in a
structure of size $\aleph_1$ iff: there exists
  an uncountable subset of that structure's universe
 such that for any two distinct $x$ and $y$ in the set,
  $\gvp(x,y)$ holds. (See Magidor and Malitz [1977].)
 
  Observe that if only $CH$ is assumed in the ground model, but not
  $\diamond_{\gw_{1}}$, the theorem would still be aplicable since
  $\diamond_{\gw_{1}}$ can be obtained in such a case by a forcing which adds
  no reals and is of size $\aleph_1$. (See Jech [1978], Exercise 22.12.)
 
  Let us see why Theorem A implies a $\Sigma^2$
  well-order of $\R$ in the generic extension. Since $CH$ is assumed, it is
  possible to find $\SP\subseteq \gw_1$ which encodes in a  natural way a
  well-order of
 $\R$ of type $\gw_1$. For example, set $\SP\subseteq \gw_1$ in such
  a way that $\SP\cap[\gw\ga,\gw\ga+\gw)$, the intersection of $\SP$ with the
  $\ga$th $\gw$-block of $\gw_1$, ``is" a subset, $r_\ga$, of $\gw$ and so that
  $\langle r_\ga$: $\ga<\gw_1\rangle$, 
is an enumeration of $\R$. Now use Theorem A to find
  a formula $\psi$ and a model $M$ of $\psi$ (with $P^M=\SP$) and a generic
  extension in which $M$ is the unique model of $\psi$. In this generic
  extension, the relation $r_\ga<r_\gb$, $\ga<\gb$, can be defined by:
\medskip
 
  {\em 
There is a model $K$ of $\psi$ where $r_\ga$ appears in $P^K$ before $r_\gb$
  does.}
\medskip

  Now, for any formula $\gvp$ in the Magidor-Malitz logic, the statement:
  ``there is a model $K$ of $\gvp$" is (equivalent to) a $\Sigma^2_2$ statement
  (see below), and hence the well-order of $\R$
  is $\Delta_2^2$ in the generic extension. (Since any $\Sigma^2_2$
linear order must be $\Delta^2_2$.)

We can start with any relation $\SP \subset \gw_1$ (not necessarily
a well-order of the reals) and get by Theorem A a generic extension $V^Q$
in which this relation is $\Delta^2_2$.

  To see the above remark, for any Magidor-Malitz formula $\gvp$,
 we encode the existence of $K\models \gvp$ as a $\Sigma^2_2$ statement
concerning $H(\gw_1)$ thus:\\
{\it There is a relation $\IN$ on $H(\gw_1)$ and a truth function
which defines a model $(H,\IN)$ of enough set theory, in which $\gw_1^H$ is
(isomorphic to)
  the real $\gw_1$, and inside $H$
 there is a model $K$ for the formula $\gvp$,
  such that:
  For any subformula $\gd(u,v)$ of $\gvp$ with parameters in $K$,
 if $X\subseteq \gw_1$ is such that for
 any
  two distinct 
$a,b\in X$ $\gd(a,b)$ holds, then there is such an $X$ in $H$ as well.}

  Now this statement ``is" $\Sigma^2_2$, and the model $K$ of $\gvp$ found in
$H$ is a real
  Magidor-Malitz model of $\gvp$, not only in the eyes of $H$.

  The second issue of the paper, the ``strong'' 
incompactness of the Magidor-Malitz
  logic, is an obvious  consequence of the fact that $\psi$ has no non
standard models.
 
  The proof of the Main Theorem involves a construction of an
 $\gw_1$ sequence of
  Suslin trees at the first stage (constructing the model $M$ of $\psi$), and
  then an iteration of posets which specialize given \Aron\ trees at the
  second stage (making $M$ the unique model of $\psi$). The main ingredient in
  the iteration is the definition of a new poset $\dS(\ST)$ for
specializing an \Aron\ tree $\ST$ without addition of new reals. 

For his well-known model of
  CH $\&$ {\it there are no Suslin trees} (SH), Jensen provides 
  (in $L$) a poset which iteratively specializes each of the
  \Aron\ trees. Each step in this iteration (including the limit stages) 
is in fact a Suslin tree.
 Both the square and the diamond are judiciously used to construct
this $\gw_2$-sequence of Suslin trees. Since forcing with a Suslin tree adds no new
reals, the generic extension satisfies CH. (See Devlin and Johnsbr\aa ten
[1974].)

In Shelah [1982] (Chapter 5) this result is 
obtained in the general and more flexible setting of
  proper-forcing iterations which add no reals. In particular, 
  a proper forcing which adds no reals and specializes a given \Aron\
  tree is defined there.

The poset $\dS(\ST)$ of
our paper is simpler than the one in Shelah [1982] because it
  involves no closed unbounded subsets of $\gw_1$, and so our paper could
  be profitably read by anyone who wants a (somewhat) simpler proof
  of Jensen's CH $\&$ SH.
 
  The paper is organized as follows:\\
Section ~\ref{s1} gives preliminaries and sets our notation. Section
~\ref{s2} shows how to construct sequences of Suslin trees such that, at
will, some products of the trees are Suslin while the others are special.
Section ~\ref{s3} is a preservation theorem for countable support iteration
of proper forcing: A  Suslin tree cannot suddenly lost its Suslinity at
limit stages of the iteration. Section ~\ref{s4} describes the poset which
is used to specialize \Aron\ trees. Section ~\ref{s5} shows that the
specializing posets of Section ~\ref{s4} can be iterated without adding
reals. In Section ~\ref{s6} we start with a given family of Suslin trees
and show how the iteration of the specializing posets obtains a model of
ZFC in which this given family is the family of {\em all} Suslin trees; all
other Suslin trees are killed. Sections ~\ref{s7} and ~\ref{s8}
are the heart of the paper and the reader may want to look there first to
get some motivation. In Section ~\ref{s8} a version of Theorem
A is proved first which suffices to answer Woodin's question, and 
then the remaining details (by now easy) are given to complete the
proof.

Concerning the related question for models where CH  does not hold, let us
report that:\\
\begin{enumerate}
\item Woodin obtained the following: Assume there is an inaccessible
cardinal $\gk$. Then there is a c.c.c forcing extension in which
$\gk=2^{\aleph_0}$ (is weakly inaccessible) and there is a
$\Delta^2_2$-well ordering of the reals.
\item Extending the methods of this paper, Solovay obtained a forcing poset
of size $2^{2^{\aleph_0}}$ such that the following holds in the extension:
\begin{enumerate}
\item $2^{\aleph_0}=2^{\aleph_1}=\aleph_2$,
\item MA  for $\sigma$-centered posets,
\item There is
a $\Delta^2_1$-well ordering of the reals.
\end{enumerate}
\item Motivated by this result of Solovay, Shelah
obtained the following: If $\gk$ is an inaccessible cardinal and GCH
holds on a cofinal segment of cardinals below $\gk$, then there is an
extension such that
\begin{enumerate}
\item $2^{\aleph_0}=\gk$, cardinals and cofinalities are not changed,
\item MA, 
\item There is
a $\Delta^2_1$-well ordering of the reals.
\end{enumerate}
\end{enumerate}

  Theorem (A) was obtained by Shelah during his visit to Caltech in 1985
and he would like to thank H. Woodin for
 asking this question and R. Solovay for encouraging conversations.
We also thank Solovay for some helpful suggestions which were incorporated
here.
The result of Section ~\ref{s6} (a model of ZFC with few Suslin trees) is
due to Abraham and appeared in fact in Section 4 of Abraham and Shelah
[1985]. (However, there the machinery of Jensen iteration of Suslin
trees was used, while here the approach of proper forcing is used.) The
poset $\dS(\ST)$  for specializing an \Aron\ tree $\ST$ was found by Abraham
who proved that any Suslin tree
 $\SS$ remains Suslin after the forcing, unless $\SS$ is
embeddable into $\ST$. As said above, $\dS(\ST)$ is simpler than the
corresponding poset $P$ of Shelah [1982], but the closed unbounded set
forcing involved in $P$ is still necessary in order to make two \Aron\ trees
isomorphic on a club.
\section{Preliminaries}\label{s1}
In this section we set our notations and remind the reader of some facts
concerning trees and forcings. All of these appear with more details in the
book of Jech [1978], or Todor\v{c}evi\`{c} [1984] or in the monograph Devlin
and Johansbr\aa ten [1974] which describes Jensen's results.

In saying that ``$\ST$ is a tree'' we intend that the height of $\ST$ is
$\gw_1$, each level $\ST_\ga$ is countable ($\ga<\gw_1)$, and every node has
$\aleph_0$ many (immediate) successors. We do not insist that the tree has
a unique root.

For a node $t\in \ST$ define its predecessor  branch by
$$
(\cdot,t)=\{s\in\ST \mid s <_{\ST} t\}.
$$
Usually it is required for limit $\gd$ that
$(\cdot,a)\not=(\cdot,b)$  for $a \not=b$ in $\ST$,
but we allow branches with more than one least upper bound.

For a node $a \in \ST$, $\level(a)$ is that ordinal $\ga$ such that $a \in
\ST_\ga$ (that is, the order-type of $(\cdot,a)$). We also say that $a$ is
of height $\ga$ in this case. $\ST\rest \ga$ is the tree consisting of all
nodes of height $< \ga$.

For a node $a \in \ST$, $\ST_a=\{x\in\ST \mid a \leq_{\ST} x \}$, is the
tree consisting of all extensions of $a$ in $\ST$.

A {\em branch} in a tree is a linearly ordered (usually downward closed)
subset. An {\em antichain} is a pairwise incomparable subset of the tree.
An \Aron\ tree is a tree with no uncountable branches. It is special if
there is an order preserving map $f:\ST\rarrow \Q,\ (x <_\ST y\
\mbox{implies} f(x)<f(y))$. A {\em Suslin} tree is one with no uncountable
antichain (and hence no uncountable chain as well). A Suslin tree has this
property that any cofinal branch (in an extension of the universe) is in
fact a generic branch. The reason being that for any dense open subset
$D\subset \ST$, for some $\ga,\ \ST_\ga \subset D$ 
(see Lemma 22.2 in Jech [1978]).

If $G \subset \ST$ is a branch of length $\gamma$, then for
$\ga<\gamma$, $G_\ga$ denotes $\ST_\ga \cap G$, and $G\rest \ga=G\cap(\ST\rest
\ga)$.

{\bf Product of trees:} The product $\ST^1\times\ST^2$ of two trees
consists of all pairs $\langle a_1,a_2\rangle$, where for some $\ga,\ a_i\in
\ST^i_\ga$. The pairs are ordered coordinatewise: $\langle
a_1,a_2\rangle < \langle a^{'}_1,a^\prime _2\rangle$ iff $a_i<_{\ST^i}a^\prime
_i$ for
both $i$'s. The product of a finite number of trees is similarly defined.

When $\langle \ST^\xi \mid \xi < \ga \rangle$ is a sequence of trees, and
$e=\langle \xi_1 \dots \xi_n \rangle$ is a sequence (or set) of indices,
then the product of these $n$ trees is denoted
$$
\ST^{(e)} = \bigtimes_{\xi \in e} \ST^\xi= \ST^{\xi_1} \times \ldots \times
\ST^{\xi_n}.
$$
This notation should not be confused with the one for their union:
$$
\ST^e = \bigcup \{ \ST^\xi \mid \xi \in e \}= \ST^{\xi_1} \cup \ldots \cup
\ST^{\xi_n}.
$$
The union of the trees is defined under the assumption that their domains
are pairwise disjoint. (It is to simplify this definition that we drop the
requirement for a unique root.)

 A {\em derived} tree of $\ST$ is formed by taking, for some $\ga < \gw_1$,
 $n$ distinct nodes, $a_1,\ldots,a_n$, and forming the product
$\ST_{a_1}\times \ST_{a_2} \times \cdots \ST_{a_n}$.

The product of a Suslin tree with itself is never a Suslin tree. And the
product of a special tree with any tree is again special.

In Devlin and Johansbr\aa ten [1974] Jensen constructs (using the diamond
$\diamond$) a Suslin tree such that all of its derived trees are Suslin
too. We will describe this construction in Section 3.1.

  Let \ST\ be an \Aron\ tree (of \height\ $\gw_1$). A function $f:\ST\raro\Q$ is
  a {\em specialization} (of \ST) if $x<_{\ST} y\Raro f(x)<_{\ST} f(y)$.
 When $f$ is
  a partial function, it is called a {\em partial specializing} function on \ST.
 
   $^n\ST_\gb$ denotes  the set of all
  $n$-tuples $\bar{x}=\seqn{x_0}{x_{n-1}}$ where $x_i\in\ST_\gb$ for all
  $i<n$.
 We also write $\bar{x}\in\ST_\gb$ instead of $\bar{x}\in\ ^n\ST_\gb$.
 $^n\ST=\cup\{ ^n\ST_\gb\mid \gb < \height \ST\}$. For $Y \subseteq
^n\ST,\ Y_\gamma=Y \cap ^n\ST_\gamma$. 
  If $x\in\ST_\gb$ and $\ga\leq\gb$,  then $x\lceil \ga$ denotes the unique
  $y\leq_{\ST} x$ with $y\in\ST_\ga$.  Similarly, for $\bar{x}\in\ST_\gb$,
  $\bar{x}\lceil\ga\eqdf\seqn{x_0\lceil
\ga}{x_{n-1}\lceil\ga}\in\ST_\ga$. If $\ga<\gb$, and
  $X\subseteq\ ^n\ST_\gb$, then $X\lceil\ga\eqdf\setm{\bar{x}\lceil
\ga}{\bar{x}\in X}$.

  Also, if $h:\ ^n\ST_\gb\raro\Q$ is a finite function, then for $\ga<\gb$,
 if the projection taking $x \in ^nT_\gb$ to $x \lceil\ga$ is one to one,
then $h\lceil
\ga:\ST_\ga\raro\Q$ is the function $h'$ defined by $h'(x\lceil
\ga)=h(x)$. And
  for a set, $H$, of finite functions, $H\lceil\ga=\setm{h\lceil\ga}{h\in H}$.

  We use similar notation for a branch, $B$, of $\ST$ denoting by
$B\lceil\ga$ the
  subset of $B$ consisting of those nodes of $B$ of $\height<\ga$. 

  Sometimes, we think of $\bar{x}\in\ ^n\ST$ as a set rather than a sequence.
  For example, when we say that $\bar{x}_1$ and $\bar{x}_2$ are disjoint: in
  this case we refer to the range of the sequences, of course. More often,
  $\bar{x}$ refers to the sequence $\seqn{x_0}{x_{n-1}}$. For example,
  $\bar{x}_1\leq_{\ST} \bar{x}_2$ means that
  $\length(\bar{x}_1)=\length(\bar{x_2})=n$, and for $i<n$, $x_{1i}\leq_{\ST}
  x_{2i}$. We do not demand that $x_i\neq x_j$.
 
  A set of $n$-tuples, $X\subseteq\ ^n\ST_\gb$, is said to be {\em
  dispersed} if for every finite $t\subseteq\ST_\gb$ there is an $n$-tuple in
  $X$ disjoint to $t$. The following Lemma is from Devlin and Johnsbr\aa
ten [1974] (Lemma 7 in Chapter VI):
\begin{lemma}
\label{l1.1}
If $\ST$ is an \Aron\ tree and
$X\subseteq\ ^n\ST$ is uncountable
  and downward closed ($\bar{y}\leq_{\ST}\bar{x}\in X\Raro\bar{y}\in X$), 
   then , for some $\gb<\gw_1$, there is an uncountable $Y\subseteq X$ such 
that:

  \ben
  \item For $\gb\leq \gga_0<\gga_1<\gw_1,\ Y_{\gga_{0}}=Y_{\gga_{1}}\lceil
\gga_0$.
  \item $Y_\gga=Y\cap\ ^n\ST_\gga$ is dispersed for every $\gb\leq\gga<\gw_1$
  (Equivalently, $Y_\gb$ is dispersed).
  \een
\end{lemma}
\section{Construction of Suslin trees} \label{s2}

  The diamond sequence, $\diamond$, on $\gw_1$, enables the construction of
  Suslin trees with some degree of freedom concerning their products.
 For example, the construction of
  two Suslin trees $\SA$ and $\SB$ such that the product $\SA\times\SB$ is
  special; or the construction of three Suslin trees $\SA, \SB$
  and $\SC$ such that $\SA\times\SB\times\SC$ is special, but $\SA\times\SB$,
  $\SA\times\SC$ and $\SB\times\SC$ are Suslin trees. This
  freedom is demonstrated in this section by showing that, given
 any reasonable prescribed
  requirement on which products are Suslin and which are special,
 the diamond constructs a sequence
  $\ssseqn{\SS(\zeta)}{\zeta<\gw_1}$ of trees satisfying this requirement.
  `Reasonable' here means that no subproduct of a Suslin product is required
  to be special.
 
  None of the {\em ideas} in this section is new, and we could have shortened
  our construction by referring to Devlin and Johnsbr\aa ten 
[1974], and leaving the
  details to the reader. We decided however to give a somewhat 
fuller presentation 
   in the hope that some readers will find it useful. The
  constructions are presented gradually, so that for the more complex
  constructions we can concentrate on the main ideas and claim that some of
  the technical details are as before.  In the following subsection 
we use the diamond to
  construct a Suslin tree such that all of its derived trees are Suslin as
  well. (Recall that a derived tree of $\ST$ has the form $\ST_{
  a_{1}}\times\ldots\times\ST_{a_{n}}$ where $a_1,\ldots,a_n\in\ST_\ga$ are distinct
  members of the $\ga$th level of $\ST$ for some $\ga<\gw_1$.)  Then we show
  the construction of two Suslin trees $\SA$ and $\SB$ such that
  $\SA\times\SB$ is special; and the last subsection gives the desired
general  construction.
 
  For the rest of this section we assume a `diamond' sequence
  $\ssseqn{S_\xi}{\xi<\gw_i}$ where $S_\xi\subseteq\xi$ and for every
  $X\subseteq\gw_1,\ \setm{\xi}{X\cap\xi=S_\xi}$ is stationary in $\gw_1$.
 
  \subsection{A Suslin tree with all derived trees Suslin} \label{s2.1}
 
  Let us recall 
the construction of a Suslin tree $\SA_\ga$. The $\ga$th level of the
 Suslin tree $\SA$, is defined by
  induction on $\ga$. In order to be able to apply the diamond to $\SA$ we
  wish to see $\SA$'s universe as $\gw_1$ and assume that  the subtree $\SA\rest\ga$ consists of
  the set of ordinals $\gw\ga$. For $\ga<\gb$ we shall require that the tree
  $\SA\rest\gb$ is an end-extension of $\SA\rest\ga$ (the reader is asked to
  forgive us for using the notation $\SA\rest\gb$ even though the tree
$\SA$ itself
has not yet been constructed).
 
  At successor stages, the passage from $\SA\rest(\mu+1)$ to $\SA\rest(\mu+2)$,
  that is, the construction of $\SA_{\mu+1}$, requires no special care: only
  that each node in $\SA_\mu$ has countably many extensions in $\SA_{\mu+1}$.
 
  At limit stage, $\gd<\gw_1$, first set $\SA\rest\gd=\bigcup_{\ga<\gd}\
  \SA\rest\ga$, and then the $\gd$th level $\SA_\gd$ is obtained by defining
(as follows) a
  countable set of branches, $\setm{b_i}{i\in\gw}$, and putting one point in
  $\SA_\gd$ above each $b_i$. The branch $b_i$ is cofinal in $\SA\rest\gd$,
  and each node in $\SA\rest\gd$ is contained in at least one $b_i$.
 
  If we only wish to construct a Suslin tree, then the diamond set
  $S_\gd\subseteq\SA\rest\gd$ is used as usual: Each node $a$ in $\SA\rest\gd$
  is first extended to some $x_0\geq a$ in $S_\gd$ (if possible), and then, in
  $\gw$ steps, an increasing sequence, $x_0<x_1<\ldots,$ is defined so that
  $\level(x_i),\ i<\gw_1$, is cofinal in $\delta$. This sequence defines one of
  our countably many branches. We see now the need for the following
  statement to hold at every stage $\delta <\gw_1$.
  $$
\mbox{\em For any }\zeta_0<\zeta_1<\delta\  
\mbox{\em  and }a\in\SA_{\zeta_{0}},\mbox{\em there is some }
  b\in\SA_{\zeta_{1}}\mbox{\em extending }a.$$
  Now let us require a little more of $\SA$ and ask that each of its
  derived trees is Suslin too (Devlin and Johansbr\aa ten [1974]).
 This variation is manifest in the construction of $\SA_\gd$ for limit
$\gd$, and is perhaps
  better described by means of a 
generic filter over a countable structure as follows.
 
  Let
  $\dP=\dP(\SA\rest\gd)$ be the poset defined thus:
  $$
  \dP=\setm{\bar{a}}{\mbox{\it
 for some }n,\ \bar{a}=\trpl{a_0}{\ldots}{a_{n-1}} 
\mbox{\it and,
  for some }\ga<\gd,\\ \mbox{\it for all }0\leq i<n,\ a_i\in\SA_\ga}
  $$
 
  The $\level$ of $\bar{a}\in\dP$ is the $\ga$ such that $a_i\in\SA_\ga$, for
  all $i<\length(\bar{a})$. A partial-order, $\bar{b}$ \extends\ $\bar{a}$, is
  defined:\\
  $\bar{b}\  \extends\  \bar{a} $ iff\\
 \level\ $(\bar{a})\leq\level(\bar{b})$,
  $\length(\bar{a})\leq\length(\bar{b})$, and $\A i<\length(\bar{a})$,
  $a_i\leq b_i$ (in $\SA\rest\gd$).\\
It is not required for $\bar{a}\in \dP$ to be one-to-one: $a_i=a_j$ is
possible, although by genericity, they will split at some stage.
 
  Now let $M$ be a countable model (of a sufficient portion of set-theory)
  which includes $\dP$ and $\SA\rest\gd$ and the diamond set $S_\gd$; and let $G$ be a
  $\dP$-generic filter over $M$. Using suitably defined dense sets, it is easy
  to see that for each fixed $i<\gw$, $b_i=$\setn{$x$}{for some $\bar{a}\in
  G,\ x=a_i$} is a branch in $\SA\rest\gd$ going all the way up to $\gd$.
  The collection $\setm{b_i}{i<\gw}$ determines $\SA_\gd$ and this ends the
  definition of $\SA$.
 
  Let $\ST=\SA_{t_{0}}\times\ldots\times\SA_{t_{n-1}}$ be any derived tree of
  $\SA$, we will prove that $\ST$ is Suslin. Let $\ga$ be the level of
  $\seqn{t_0}{t_{n-1}}$ in $\SA$ (so $t_i \in \SA_\ga$ for all $i < n$).
 Let $E\subseteq\ST$ be any dense open subset.
  By the diamond property, using some natural encoding of $n$-tuples of
  ordinals as ordinals, for some limit $\gd>\ga$,
  $E\cap(\ST\rest\gd)=E\cap\gd=S_\gd$, and $E\cap \gd$ is dense open in
  $\ST\rest\gd$. We must prove that every $\bar{x}\in\ST_\gd$ is in $E$, in
  order to be able to prove that an arbitrary antichain in $\ST$ is countable.
   $\bar{x} \in \ST_\delta$
 has the form $\bar{x}=\seqn{x_0}{x_{n-1}}$ where $x_k\in
  \SA_{t_{k}}$. Since the $t_k$'s are distinct, the $x_k$'s give distinct
  branches of $\SA\rest\gd$. Recall the $\dP(\SA\rest\gd)$ generic filter $G$
  (over $M$) used to define $\SA_\gd$, and let $b_{i(k)}$ be the branch of
  $\SA\rest\gd$ which gave $x_k$. If for some
$\bar{a}=\seqn{a_0}{a_{l-1}} \in G$ 
the
  $n$-tuple $\seqn{a_{i(0)}}{a_{i(n-1)}}$
is in the dense set $E$, then
  $\bar{x}$ which is above this $n$-tuple must be in $E$ too.
 
  The existence of such $\bar{a}$ in $G$ is a consequence of the following
  density argument:  Let $D$ contains all those $\bar{a}\in\dP(\SA\rest\gd)$
  for which (1) $i(k)<\length(\bar{a})$ for every $k<n$, and either the
  subsequence $s=\seqn{a_{i(0)}}{a_{i(k)}}$ of $\bar{a}$ is not in $\ST$,
 or else
  $s\in E$. $D$ is dense in $\dP$ and $D\in M$, because $S_\gd$ is in
$M$. So that $D\cap G\neq 0$, and
  any $\bar{a}\in D\cap G$ of height  $>\ga$ is as required.
 
  \subsection{The case of two trees} \label{s2.2}
 
  The construction in the previous section is combined now with the
  construction of a special \Aron\  tree to yield two Suslin trees $\SA$ and
  $\SB$ such that:
  \ben
  \item Each derived tree of $\SA$ and of $\SB$ is Suslin,
  \item $\SA\times\SB$ is a special tree.
  \een
 
  We commence by recalling the construction of a special 
\Aron\  tree $\SA$ together with a
  strictly increasing $f:\SA\raro \Q$. In the inductive definition, $\SA\rest a$
  and $f\rest\ga= f\rest(\SA\rest \ga)$, $\ga < \gw_1$,
 are defined so that the following
  hold:
$$
  ({\it 1})\ \mbox{\em For any }a\in\SA\rest\ga\ \mbox{\em and rational }\gep>0,
\mbox{\em and ordinal }\tau<\ga\ \mbox{\em such
  that }\height(a)<\tau,$$ $$\mbox{\em there is an extension }b>a \
\mbox{\em in }\SA_\tau\ \mbox{\em such
  that }0< f(b)-f(a)<\gep.$$

  This condition is needed at a limit stage $\ga$, if we don't want our cofinal
  branches to run out of rational numbers; it enables the assignment of
  $f(e)\in\Q$ for $e\in\SA_\ga$, but it requires some care to keep it true at
  all stages.
 
  There is nothing very special at successor stages: Since we assume that
each node has $\aleph_0$ many successors, 
condition (1) above may be achieved by assigning to these successors of $e$
all the possible values of rational numbers  $>f(e)$ (a forcing-like
description of the successor stage is also possible---see below).

  For a limit $\ga < \gw_1$, it seems again convenient to formulate the
  construction of $\SA_\ga$ and $f\rest\SA_\ga$, in terms
  of a generic filter $G$ over a countable structure $M$. So given
$A\rest\ga$ and $f\rest\ga$, a countable poset
  $\dq=\dq(\SA\rest\ga,f\rest\ga)$ is defined first.
 
  \begin{defn}  Let $\dq$ be the collection of all pairs $(\bar{a},\bar{q})$
  of the form $\bar{a}=\seqn{a_0}{a_{n-1}}$  $\bar{q}=\seqn{q(0)}{q(n-1)}$
  such that:
  \ben
  \item $\bar{a}$ is an $n$-tuple in $\SA_\xi$ for some $\xi<\ga$.
  \item $q(i)\in\Q$ and $f(a_i)<{q}(i)$ for all $i<n$.
  \een
\end{defn}
 
  Intuitively, ${q}(i)$ is going to be the value of $f(b)$ for $b\in
  \SA_\ga$ defined by the `generic' branch $\setm{a_i}{(\bar{a},\bar{q})\in
  G}$. The order relation on $\dq$ is accordingly defined:
  $(\bar{a}_2,\bar{q}_2)\extends(\bar{a}_1,\bar{q}_1)$ iff
  $\bar{a}_2\ \extends\ \bar{a}_1$  and $\bar{q}_1$ is an
  initial sequence of $\bar{q}_2$ (that is,
  $\length(\bar{q}_1)\leq\length(\bar{q}_2)$, and for $i<\length(\bar{q}_1),\
  \bar{q}_1(i)=\bar{q}_2(i)$).

  Now we turn to the construction of two Suslin trees $\SA$ and $\SB$ such
  that all the
 derived trees of $\SA$ and $\SB$ are Suslin and yet $\SA\times\SB$
  is special. In this case both $\SA\rest\ga$, $\SB\rest\ga$ and the
  specializing function $f:(\SA\rest\ga)\times(\SB\rest\ga)\raro \Q$ are
  simultaneously constructed by induction. There are three jobs to do at the
limit
  $\ga$th stage: (i) to ensure that the cofinal branches of $\SA\rest\ga$ and
  of its derived trees all pass through the diamond set $S_\ga$. (ii) To ensure
  the similar requirement for $\SB\rest\ga$. (iii) To specialize
  $(\SA\rest\ga+1)\times(\SB\rest\ga+1)$. It turns out that it suffices to
  take care of (iii) in a natural way - and genericity will take care of
the two other 
  requirements, thereby ensuring that
 $\SA$ and $\SB$ and their derived
  trees are Suslin.

The inductive requirement ({\it 1}) is needed here too, but in fact an even
stronger requirement will be used:

 {\it (2)} If $\bar{a},\ \bar{b}$ are $n$ and $m$ tuples in
$\SA\rest\ga$ and $\SB\rest\ga$, and for $\gt < \ga$
 $\bar{c}$ is an $n$-tuple in $\SA_\gt$
extending $\bar{a}$, and if $q:\ n \times m \raro \Q$ is such that
$\forall i,j\ f(a_i,b_j) < q(i,j)$, THEN
there is an $n$-tuple $\bar{d}$ in $\SB_\gt$, extending $\bar{b}$ 
such that
$\forall i,j\ f(c_i, d_j)<q(i,j)$. (A similar requirement is made for
$\bar{c}$ in $\SB_\gt$.)
 
In fact, by taking smaller $q(i,j)$, it even follows that
 for any finite $D \subset B_\gt$
there is an $n$-tuple $\bar{d}$ in $\SB_\gt$, extending $\bar{b}$ and
disjoint to $D$, such that
$\forall i,j\ f(c_i, d_j)=q(i,j)$. (A similar requirement is made for
$\bar{c}$ in $\SB_\gt$.)
  Again, we only describe the limit case, and leave the details of the
  successor case to the reader (take care of condition ({\it 2})).
 So assume $\ga<\gw_1$ is a limit ordinal and
  $\SA\rest\ga\ (=\cup_{\mu<\ga}\SA\rest\mu),\ \SB\rest\ga$ and $f\rest\ga$ are
  given. 

 Let $\R=\R((\SA\rest\ga)\times(\SB\rest\ga), f\rest\ga)$ be the poset
  defined in the following: $(\bar{a},\bar{b},\bar{q})\in\R$ iff for some
  $\mu<\ga$, $\bar{a}$ is an $n$-tuple in $\SA_\mu$, $\bar{b}$ is an
  $m$-tuple in $\SB_\mu$, and $\bar{q}:n\times m\raro\Q$, are such
  that for all $0\leq i<n,\ 0\leq j<m$, $f(a_i,b_j)<\bar{q}(i,j)$. Extension
  is defined naturally.
 
  If $M$ is now a chosen countable structure containing all the above, and the
  diamond $S_\ga$ in particular, then pick an $\R$-generic filter, $G$, over
  $M$ and define the $\ga$th levels $\SA_\ga$ and $\SB_\ga$ and extend $f$ on
  $\SA_\ga\times\SB_\ga$ in the following way: For each $i$,
  $$
  \begin{array}{l}
  u_i=\mbox{\it a node above the branch\ }
\setm{x}{\mbox{\it For some }(\bar{a},\bar{b},\bar{q})\in G,\ \ a_i=x}\\
  v_i=\mbox{\it a node above the branch\ }
\setm{y}{\mbox{\it For some }(\bar{a},\bar{b},\bar{q})\in G,\ \ b_i=y}\\
  f(u_i,v_j)=\bar{q}(i,j),\mbox{\it where } (\bar{a},\bar{b},\bar{q})\in G\mbox{
  \it for some }\bar{a},\bar{b}
  \end{array}
  $$
 
  By condition ({\it 2}), it is clear that this local forcing
 adds branches above every node,
and that condition ({\it 2}) continues to hold for $\ga+1$. So
 the product of the trees $\SA$ and $\SB$ thus 
obtained is
  special; but why are $\SA$ and $\SB$ and each of their derived trees
  Suslin? To see that, we argue that if we restrict our attention to $\SA$,
  for example, then it is in fact the construction of a Suslin tree given in
  subsection \ref{s2.1} which describes $\SA$. For this aim we will define for
  each limit $\ga$ a projection $\Pi_{\SA}$ from the poset $\R$ used in the
  construction of $(\SA\rest\ga+1)\times(\SB\rest\ga+1)$ onto the poset
  $\dP(\SA\rest\ga)$ used in \ref{s2.1}. Simply set
  $\Pi_{\SA}(\bar{a},\bar{b},\bar{q})=\bar{a}$. We must check the following
  properties which ensure that the projection of the $\R$-generic filter is
a $\dP(\SA\rest\ga)$-generic filter:
  \ben
  \item $\Pi_{\SA}$ is order preserving: $x_0\leq_{\mbox{\small \R}}
  x_1\Raro\Pi_{\SA}(x_0)\leq_{\dP}\Pi_{\SA}(x_1)$.
  \item Whenever $p\in\dP$ extends $\Pi_{\SA}(x_0)$, there is an extension $x_1$
  of $x_0$ in $\R$ such that $\Pi_{\SA}(x_1)=p$.
  \een
 
  This is not difficult to prove by ({\it 2}).
 
  \subsection{$\gw_1$ many trees} \label{s2.3}
 
  In this section (waving our hands even harder) we extend the previous
  construction to $\gw_1$ many Suslin trees with any reasonable requirement on
  which trees are Suslin and which are special.
 
  \begin{theorem}
 
  Assume $\diamond_{\gw_{1}}$. Let \spa\ (for special) be a collection of
  non-empty
finite subsets of $\gw_1$ closed under supersets,
 and let \sua\ be those non-empty 
finite sets $e\subset
  \gw_1$ which are not in $\spa$.
 Then there is a sequence of $\gw_1$-trees
  $\sseqn{\SA^\zeta }{\zeta<\gw_1}$ such that for a finite set
$e=\fsetn{\zeta_1}{\zeta_n}$
  \ben
  \item If $e\in\spa$, 
  $\SA^{(e)}\stackrel{Def}{=}
\nseqn{\SA^{\zeta_1}\times}{\times\SA^{\zeta_n}}$ is special.
  \item If $e\in\sua$, $\SA^e \stackrel{Def}{=}
\SA^{\zeta_1}\cup\ldots\cup\SA^{\zeta_n}$ and
  all of its derived trees are Suslin.
 \een
  \end{theorem}
 
  \proof By induction on $\ga<\gw_1$, the sequence
  $\setm{\SA^\zeta\rest\ga+1}{\zeta<\ga}$ is defined together with specializing
  functions $f_e\rest\ga+1$ for $e\in\spa,\ e\subseteq\ga$. $f_e$ is, of course,
  a specializing function from $\SA^{(e)}$ into $\Q$.
It is convenient to require that $f_e\rest\ga+1$ is only defined on the
$\gb$ levels of the product tree for $\gb > \max (e)$.
 
The definition of the trees requires some notations and preliminary
definitions. Let $\ga < \omega_1$ be any ordinal---successor or limit, and
assume that
  $(\SA^\zeta\rest\ga)$ for $\zeta<\ga$, and $f_e\rest\ga$ for $e\subseteq\ga$
  in \spa\ are given.

Let us define, for any finite $d \subseteq \ga$,
$$
\dP(\SA^d\rest\ga)=\bigtimes_{\xi \in d} \dP(\SA^\xi\rest\ga)
$$
in the following. $a=\langle\ol{a}^\xi \mid \xi \in d \rangle \in
\dP(\SA^d\rest\ga)$ if for some $\mu < \ga$ for all $\xi \in d,\ 
\ol{a}^\xi$ is an $n_\xi$-tuple in $\SA^\xi_\mu$ enumerated as follows:
$\ol{a}^\xi=\langle a^\xi_i \mid i \in I_\xi \rangle$ where $ I_\xi =
I_\xi(a) \subset \omega$ is a finite set of size $n_\xi$.

Extension is naturally defined in $\dP(\SA^d\rest\ga)$: $b$ extends $a$ iff
for all $\xi \in d$,\ $I_\xi(a) \subseteq I_\xi(b)$ and for every $i \in
I_\xi(a),\ a^\xi_i < b^\xi _i$ in $\SA^\xi \rest \ga$.

We need one more definition. For $a \in \dP(\SA^d \rest \ga)$, and $e
\subseteq d$ with $e \in \spa$, let us say that $q^e$ {\em bounds} $a$
iff  $q^e$ is a function
 $$
q^e:\ I^{(e)}=\bigtimes_{\xi \in e} I_\xi(a) \rarrow \Q
$$
such that for every $\ol{i}=\langle i_\xi \mid \xi \in e \rangle \in
I^{(e)},$
\[ f_e(\langle a^\xi_{i_\xi} \mid \xi \in e \rangle )
< q^e(\ol{i}).\]

Now we can formulate the inductive property of the trees and functions at
the $\ga $ stage:\\
$
(3_\alpha)$  {\it If } $a \in \dP(\SA^d \rest \ga),$ {\it where}
$ d \subseteq \ga$ {\it is finite, 
 and if} $\langle q^s \mid s \subseteq d,\ s \in \spa \rangle$ 
{\it is
such  that each} $q^s$ {\it bounds} $a$, {\it then for every }
 $e \subseteq d$ {\it
with } $e \in \sua,$ {\it and } $b$ {\it extending } $a\rest e$ 
 {\it in }
 $\dP(\SA^e \rest \ga),$ {\it there is} $b_1>a$ {\it in}
 $\dP(\SA^d \rest \alpha)$ {\it such that each}
 $q^s$ {\it bounds } $b_1,$ {\it and } $b_1 \rest e =b.$
\\
Observe that this property makes sense for $\alpha$'s which are  limit as
well as for  $\alpha$'s which are successor ordinals.

Let us now return to the inductive definition of the trees.\\
{\bf Case 1} $\ga$ is a limit ordinal. In this case we first take the union
of the trees and functions obtained so far.
So for each $\xi<\ga$, and $e \subseteq \ga$ in  $\spa$:
$$
\SA^\xi \rest \ga = \bigcup_{\mu < \ga} \SA ^\xi \rest \mu,\ {\it and }
f_e \rest \ga = \bigcup _{\mu < \ga}  f_e\rest \mu,\ 
$$
Then we add the $\ga$-th levels and extend $f_e$ according to the following
procedure.

A countable poset $\R=\R_\ga$ is defined as a convenient
  way to express how the $\ga$-branches are added to each
  $(\SA^\zeta\rest\ga)$, enabling the definition of $\SA^\zeta_\ga$ and of the
  extensions of the $f_e$'s.

 A condition in $\R$ gives finite information on
  the branches and the values of the appropriate $f_e$'s. So a condition
 $r\in\R$ has two components:
$r=\{a,\ol{q}\}$, where for some $\mu=\mu(a)<\ga$, $a$ gives information on
the intersection of the (locally) 
generic $\ga$-branches with the $\mu$ level, and
$\ol{q}$ tells us the future rational values on products of 
these branches. Formally, we require
that for some finite $d=d(r) \subseteq \ga$
$$
\begin{array}{ll}
(1) & a \in \dP(\SA^d\rest\ga),\\
(2) & \ol{q}=\langle q^e \mid e \subseteq d,\ e \in \spa \rangle,\ \mbox{and
for each } e \subseteq d\ \mbox{in } \spa,\ q^e \ \mbox{bounds } a.
\end{array}
$$
  In plain words, $r\in\R_\ga$ has a finite domain $d\subseteq\ga$ on which it
  speaks. For $\zeta\in d$, $a^\zeta_i$ is the intersection with
$\SA^\zeta_{\mu}$, of the proposed
  $i$th branch added to $\SA^\zeta\rest\ga$, and 
  $q^e$ gives information on how to specialize those trees required to be
  special. So if $e\subseteq d$ is in \spa, then $q^e$ gives rational
  upper bounds to the range of the specializing function $f_e$ on the branches
  added to $\SA^{(e)}\rest\ga$.
 
  We write $d=d(r),\ \bar{a}=\bar{a}(r),\ \bar{q}=\bar{q}(r),\ \mu=\mu(r)$ etc.
 to denote the components of $r\in\R$.
 
   A countable $M$ is chosen with $\R_\ga$, $S_\ga$, the trees so far
constructed and so on in $M$, and an $\R_\ga$ generic filter $G$ over $M$ is
  used to define the branches and the new values of $f_e$.
  \ben
  \item For $\zeta<\ga$ and $i<\gw$, $u_i^\zeta=\setm{x}{\mbox{\it For some }
r\in
  G,\ a_i^\zeta=x}$, is the $i$th branch added to $\SA^\zeta\rest\ga$. This
  determines $\SA_\ga^\zeta$.
  \item For $e\in \spa$, we define $f_e$ on the $\ga$th level of
  $\SA^{(e)}$, as follows.
  Any $\ga$-level node, $w$, of $\SA^{(e)}$ has the form
  $\seqn{u_{i_1}^{\zeta_1}}{u_{i_n}^{\zeta_n}}$ where
  $e=\fsetn{\zeta_1}{\zeta_n}$; then $f_e(w)=q^e(i_1,\ldots,i_n)$,
  where $q^e$ comes from $G$,
  (That is $q^e=\bar{q}(r)^e$ for some $r\in G$).
\een
 
  As evidenced by $f_e$, $\SA^{(e)}$
 becomes a special tree for $e\in \spa$. When
  $e\not\in\spa$,  $\SA^e$ is Suslin and so are all of its derived trees.
  It is here that the assumption $e\not\in\spa\Raro\mbox{ for } e'\subseteq e,\
  e'\not\in \spa$ is used. We must prove that for $e\not\in\spa$, the
  construction of $\SA^e$ follows the specification described in subsection
  \ref{s2.1}. To do that, observe that for $e\not\in\spa$ the map $\Pi$ taking
  $r$ to $\ssseqn{\bar{a}^\zeta(r)}{\zeta\in e}$ is a projection of $\R$ onto
  $\dP(\SA^e\rest\ga)$.

{\bf Case 2} is for $\ga$ a succesor ordinal. Put $\ga = \rho+1$. Not only
the $\ga$th level has to be defined for all existing trees, but a new tree,
$\SA^\rho$, and new functions must be introduced. The definition of the new
functions $f_e$ with $\rho \in e$ is somewhat facilitated by our assumption
that these are only defined on the $\ga$th level.
  \section{Suslin-tree preservation by proper forcing} \label{s3}
 
  In this section it is shown that any Suslin tree $\SS$ remains Suslin 
  in a countable support iteration, if each
  single step of the iteration does not destroy the Suslin property of $\SS$.
We assume our posets are separative: If $q$ does not extend $p$ then some
extension of $q$ is incompatible with $p$.
 
  Let $\vec{\dP}=\ssseqn{\dP_\ga}{\ga\leq\gb}$ be a countable support iteration
  of length $\gb$ (limit ordinal) of proper forcing posets; where
  $\dP_{\ga+1}=\dP_\ga*\dq_\ga$ is a two step iteration: $\dP_\ga$ followed by
  $\dq_\ga$. The following preservation theorem holds for $\vec{\dP}$.
 
  \begin{theorem} \label{t4.1}
  Let $\SS$ be a Suslin tree of \height\ $\gw_1$; suppose for every $\gb'<\gb$,
  $\SS$ remains Suslin in $V^{\dP_{\gb'}}$. Then \SS\ remains Suslin in
  $V^{\dP_{\kern -2pt \gb}}$ as well.
 \end{theorem}
 
  \proof To show that every antichain of \SS\ is countable, it is enough to
  prove that:\\

  { \em For any dense open set $E\subseteq \SS$ there is a level $\SS_\gl$,
  $\gl<\gw_1$,  such that $\SS_\gl\subseteq E$}.
 
  So let $\bE$  be a $\dP_\gb$ name of a dense open subset of \SS.  Fix some
  countable elementary submodel, $M\prec H(\gk)$, such that \SS, $\gb,\
  \dP_\gb,\ \bE$ etc. are in $M$, where $\gk$ is ``big enough". ($H(\gk)$ is the
  collection of all sets of cardinality hereditarily $<\gk$. In fact, all that
  is needed is that $M$ reflects enough of $V$ to enable the following
  constructions and arguments to be carries out.)
  \bd
  \item{$\bullet$} Let $\gl=M\cap\gw_1$. $\gl$ is a countable ordinal.
  \item{$\bullet$} Let $\ssseqn{\gb(i)}{i<\gw}$ be an increasing
  $\gw$-sequence of ordinals in $\gb\cap M$ and cofinal in $\gb\cap M$.
  \item{$\bullet$} Let $\setm{b_n}{n\in\gw}$ be an enumeration of $\SS_\gl$.
  \ed
 
  We will produce a condition $q\in\dP_\gb$ (extending some given condition)
  such that for every $n\in\gw$,
  $q\force b_n\in \bE,$ and thus
  $$
  q\force \SS_\gl\subseteq \bE
  $$
 
  First, a sequence $q_n\in\dP_{\gb(n)}$, and $\bp_n$ is inductively
  constructed such that the following holds:
  \ben
  \item $q_n$ is an $M$-generic condition for $\dP_{\gb(n)}$ and
  $q_{n+1}\rest\gb(n)=q_n$.
  \item $\bp_n$ is  a \name\ in $V^{\dP_{\kern -2pt\gb(n)}}$, forced to be a condition in
  $\dP_\gb\cap M$.
  \item
  \bd
  \item{(a)} $q_n\Pbnforce\ \bp_n\rest\gb(n)$ {\em is in the canonical generic
  filter}, $G_n$.
  \item{(b)} $q_n{\Pbnforce}\ \bp_n$ \extends\ $\bp_{n-1}$ {\em in} $\dP_\gb$.
  \item{(c)} $q_n\Pbnforce$ ($\bp_n\Pbforce\ b_n$ {\em is above some member of }
  $\bE$).
  \ed
  \een
 
  (Recall that the canonical generic filter $G$ is  defined so that $q\force
  q\in G$ for every $q$).
 
  Suppose for a moment that we do have such sequences, and this is how $q$ is
  obtained: $q=\cup_{n<\gw}q_n$. Then $q\in\dP_\gb$ extends each $q_n$, since
  $q_{n+1}\rest \beta(n)=q_n$ is assumed for all $n$.
We also have the following:
 
  \begin{claim} $q\Pbforce \SS_\gl\subseteq \bE$.
\end{claim}
 
  \proof This is a consequence of 1-3 obtained as follows. Let $\bG$ be the name
  of the $\dP_\gb$ canonical generic filter. It is enough to show that for
  each $n$ 
  $$
  \leqno(*)\hspace{1cm} q\Pbforce\ \bp_n\in \bG,
  $$
  because then we use 3(c) to deduce that $q$ forces $b_n$ to be in $\bE$.
  To prove (*) we observe that
  \ben
  \item $q\Pbforce\ \bp_n$ {\em extends} $\bp_m$ {\em in } $\dP_\gb$ {\em for }
  $m<n$, and
  \item $q\Pbforce\ (\bp_m\rest\gb(m))\in \bG_m$ {\em for all } $m<\gw$.
  \een
Hence:
  
 \[ q\Pbforce (\bp_m\rest\gb(n))\in \bG_n \mbox{\em for all } m<n.\]
 
  From this it follows that for any $q'$ extending $q$ in $\dP_\gb$, if $q'$
  determines $\bp_m$, that is for some $p\in\dP_\gb$, $q'\force \bp_m=p$, then
  $q'\force p\rest\gb(n)\in \bG_n$, and hence $q'$ extends $p\rest\gb(n)$, for
  all $n$'s, and thus $q'$ extends $p$.  Thus $q'\Pbforce\
 \bp_m\in \bG$. This is
  so for an arbitrary extension of $q$ which determines $\bp_m$, and hence (*):
  $q\Pbforce\ \bp_m\in \bG$,
 
  Return now to the construction of the sequences. Suppose $\bp_n$ and $q_n$
  are constructed (or that we are about to start the construction).
 
  In order to describe $\bp_{n+1}$ and $q_{n+1}$ (in that order), imagine a
  generic extension $V[G_n]$ of our universe $V$, where
  $G_n\subseteq\dP_{\gb(n)}$ is a generic filter containing $q_n$. Then
  $M[G_n]$ can be formed; it is the $G_n$-interpretation of all $\dP_{\gb(n)}$
  names in $M$. Then $M[G_n]\prec H(\gk)[G_n]$. \SS\ is still a Suslin tree in
  $V[G_n]$ by our assumption.
 
  In $V[G_n]$, $\bp_n$ is realized as a condition denoted $p_n$;
  $p_n\in\dP_\gb\cap M$, and $p_n\rest\gb(n)\in G_n$ by the inductive
  assumption in 3(a).
 
  Since $\bE$ is forced to be dense in \SS, for any $s\in\SS$ and $p\in\dP_\gb$,
  there are $s\leq_{\SS} s'$ and an extension $p'$ of $p$ in $\dP_\gb$ such
  that
  $$
  \leqno(**)\hspace{1cm} p'\Pbforce\ s'\in \bE
  $$
 
  Moreover, by genericity of $G_n$, we may require that
  $$
  p\rest\gb(n)\in G_n\Raro p'\rest\gb(n)\in G_n.
  $$
 
  Thus, the set $F$ of $s'\in \SS$ for which there is $p'\in\dP_\gb$ extending
  $p_n$ with $p'\rest\gb(n)\in G_n$ and satisfying 
($\ast \ast$) is dense in \SS\ and is
  (defined) in $M[G_n]$.
 
  Now \SS\ is a Suslin tree in $M[G_n]$, and hence every branch of
  $\SS\rest\gl$ of length $\gl$ is $M[G_n]$ generic. (Recall $\gl=\gw_1\cap
  M$).  Thus $b_{n+1}$ (the ($n+1$)th node of $\SS_\gl$) is above some node in
  $F$; and it is possible to pick $p_{n+1}$ in $\dP_\gb\cap M$ extending $p_n$ with
  $p_{n+1}\rest\gb(n)\in G_n$ and such that $p_{n+1}\Pbforce\ b_{n+1}\in \bE$.
 
  This description of $p_{n+1}$ made use of the $\dP_{\gb(n)}$-generic filter
  $G_n$. Back in $V$, we define $\bp_{n+1}$ to be the {\em name} of that
  $p_{n+1}$ in $V^{\dP_{\gb(n)}}$ (and so, evidently, in $V^{\dP_{\gb(n+1)}}$).
 
  Next we define $q_{n+1}$. We demand the following from $q_{n+1}$:
  \ben
  \item $q_{n+1}$ is an $M$-generic condition for $\dP_{\gb(n+1)}$, and
  $q_{n+1}\rest\gb(n)=q_n$.
  \item $q_{n+1}\force_{\dP_{\gb(n+1)}}\ \bp_{n+1}\rest\gb(n+1)\in G_{n+1}$.
  \een
 
  The existence of $q_{n+1}$ satisfying (1) and (2) is a general fact about
  proper forcing. It is a consequence of the following statement, which can
  be proved by induction on $\gb_2$:
 
  Suppose $\gb_1<\gb_2\leq\gb$, and $q_1$ is an $M$-generic condition over
  $\dP_{\gb_{1}}$, and $\bp$ is a name in $V^{\dP_{\gb_1}}$ such that
  $q_1\force_{\dP_{\gb_1}}$ $\bp\in\dP_{\gb_2}\cap M$ {\em and} $\bp\rest\gb_1$
  {\it is in the canonical $\dP_{\gb_1}$ generic filter.}
  Then there is an $M$-generic condition over $\dP_{\gb_2}$, $q_2$, such that
  $q_2\rest\gb_1=q_1$, and $q_2\Pbforces{2}\ \bp\rest \gb_2$
 {\it is in the canonical
  $\dP_{\gb_{2}}$ generic filter.}
  \section{ How to specialize \Aron\ trees without adding 
reals} \label{s4}
 
  The forcing notions which turn a given \Aron\ tree into a special tree,
  naturally fall into two categories: those which use finite conditions and
  satisfy properties such as the c.c.c., and those which use infinite
  conditions and have nice closure properties.
 
  In this section we describe how infinite conditions can be used to
  specialize an \Aron\ tree, without addition of new countable sets, and
  how to iterate such posets.

  In a moment we will define the poset $\dS(\ST)$ used to specialize an \Aron\
  tree \ST.  Meanwhile, let us see what are the problems with the direct
  approach, which takes the poset $\dS_1$ of all specializing functions $f$
  defined on some downward closed countable subtree of the form $\ST\rest
  \ga+1$. To see that this poset collapses $\gw_1$ in forcing, look at the
  following dense open sets, defined for $n<\gw$.
  $$
  D_n=\setm{f\in\dS_1 }
 {f\ 1\mbox{\it
 is defined on }
  \ST\rest\ga+1, \mbox{\it and for
 every }x\in \ST_\ga,\ f(x)\geq n}
  $$
 
  Clearly $D_n$ is dense open.
So for every $n < \gw$ there is an $\ga < \gw_1$ such that some $f \in D_n$
defined on $\ST \rest \ga+1$ is in the generic filter $G$.
   But if $\gw_1$ is {\em not}
  collapsed, the generic filter must contain a condition which is
  simultaneously in every $D_n$, and this is a contradiction.
 
  Thus, there must be some limitation on the growth of the generic
  specializing function.  We may try the following poset: $\dS_2$ consists of
  all specializing $f:\ST\rest \ga+1\raro\Q$, $\ga<\gw_1$, such that
  $$
  \A\ga_0<\ga\ \A\bar{x}\in\ST_{\ga_0}, \mbox{ if } f(\bar{x})<\bar{q}\mbox{
  then for some }\bar{y}\in\ST_\ga,\ \bar{x}<\bar{y}\mbox{ and }
  f(\bar{y})=\bar{q}.
  $$
  Here, $\bar{x}$ is an $n$-tuple of nodes,  $\bar{q}$ is an $n$-tuple of
rational numbers, and $f(\bar{x})<\bar{q}$ is a shorthand for:
  $f(x_i)<q_i$ for all $i$'s.

  If we assume that every node in \ST\ has infinitely many successors, it is
  not difficult to see that any condition in $\dS_2$ can be extended to any
  height.  Therefore forcing with $\dS_2$ specializes \ST. If \ST\ is a Suslin
  tree such that each derived tree of \ST\ is Suslin too, then $\dS_2$ adds no
  new countable sets. 
  If \ST\ is an arbitrary \Aron\ tree, however,
 $\dS_2$ may collapse $\gw_1$.  Such is
  the case when \ST\ has the form $\ST_1\cup\ST_2$, a disjoint union of $\ST_1$
  with a copy $\ST_2$ of itself.  Let $i:\ST_1\raro\ST_2$ be the map which
  takes a node in $\ST_1$ to its copy in $\ST_2$. Define
  $$
  D_n=\left\{f\in\dS_2|\begin{array}{l} dom(f)=\ST\rest\ga+1\mbox{ and }
  \A x\in\ST_{1,\ga}\\  f(x)>n\mbox{ or } f(i(x))>n\end{array}\right\}
  $$
 
  Again, $D_n$ is seen to be dense open; and since the generic function
  contains a condition in $D_n$ but there is no condition in the intersection
  of all the $D_n$'s, $\gw_1$ must be collapsed. This shows that the
  limitations one imposes on the growth of the generic specializing function
  must have a different character.
 
  We are now going to define the poset $\dS=\dS(\ST)$ used to specialize a
  given \Aron\ tree \ST. A condition $p\in\dS$ is a pair $p=(f,\gG)$ where 
$f$ is
  a countable partial specialization of \ST, called an ``approximation"; and
  $\gG$ is an uncountable object, called a ``promise". Its role is to ensure
  that $\dS$ is a proper poset. $\gG$ consists of ``requirements", so we have
  to explain what these are first. Throughout, \ST\ is a fixed \Aron\ tree.
 
  \begin{defn} \label{d4.1} (1) We say that $H$ is a {\em requirement} (of
  \height\ $\gga<\gw_1$) iff for some $n=n(H)<\gw$, $H$ is a set of finite
  functions of the form $h:\ST_\gga\raro\Q$, with $dom(h)\in\ ^n\ST_\gga$.\\
  \nin(2) An {\em approximation} (on \ST) is a partial specializing function
  $f:\ST\rest(\ga+1)\raro\Q$; that is, an order-preserving function defined on
  $\bigcup_{\zeta\leq\ga}\ST_\zeta$ into the rationals. The countable ordinal
  $\ga$ is called $\last(f)$. 

We say that a finite function $h:\ST_\ga\raro\Q$
  {\em bounds} $f$ iff $\A x\in dom(h) (f(x)<h(x))$. More generally, for
$\gb \geq \ga=\last(f),\ h: \ST_\gb \raro \Q$ bounds $f$ iff $\forall x \in
\dom(h) (f(x\lceil\ga)<h(x))$ (i.e., if $h\lceil\ga$ is defined, then 
$h\lceil\ga$ bounds $f$).\\
  \nin(3) An approximation $f$ with $\last(f)=\ga$ is said to {\em fulfill}
  requirement $H$ of \height\ $\ga$ iff for every finite $t\subseteq\ST_\ga$
  there is some $h\in H$ which bounds $f$ and such that $dom(h)$
 is disjoint to $t$.\\
  \nin(4) A {\em promise} $\gG$ (for \ST) is a function
  $\sseqn{\gG(\gga)}{\gb\leq\gga<\gw_1}$ ($\gb$ is denoted $\gb(\gG)$) such
  that
  \bd
  \item{(a)} $\gG(\gga)$ is a countable collection of requirements of
 height\ $\gga$. There is
a fixed $n$ such that $n=n(H)$ for all $\gga$ and $H \in
\gG(\gga)$.
  \item{(b)} For $\gga\geq\gb$, each $H\in\gG(\gga)$ is dispersed. That is,
  for every finite $t\subseteq\ST_\gga$ for some $h\in H$, $t\cap
  dom(h)=\emptyset$.
  \item{(c)} For every $\gb\leq\ga_0<\ga_1<\gw_1$,
  $$
  \gG(\ga_0)=\setm{(X\lceil\ga_0)}{X\in\gG(\ga_1)}
  $$
  \ed
  \nin(5) An approximation $f$ fulfills promise $\gG$ iff
  $\last(f)\geq \gb(\gG)$, and $f$ fulfills each requirement $H$ in
  $\gG(\last(f))$.
\end{defn}
 
  \begin{defn}[of $\dS(\ST)$]  
\label{d4.2}
For any \Aron\ tree \ST\ define
$\dS=\dS(\ST)$ by
  $p=(f,\gG)\in\dS$ iff $f$ is an approximation on \ST,  $\gG$ is a
  promise, and $f$ fulfills $\gG$.
 
  The partial order is naturally defined: $p_1=(f_1,\gG_1)\extends\
  p_0=(f_0,\gG_0)$ iff $f_0\subseteq f_1$ and
  $\gG_0\rest(\gw_1-\last(f_1))\subseteq\gG_1$. That is, any requirement of
  height $\gga\geq\last(f_1)$ in $\gG_0$ is also in $\gG_1$.
 
  If $p=(f,\gG)$ is a condition in $\dS$ we write $f=f(p),\ \gG=\gG(p)$. In an
  abuse of notation, we write $\last(p)$ for $\last(f(p))$, and $p(x)$ instead
  of $f(x)$. We also call $\last(p)$ `the \height' of $p$. 
(Recall that $f(p)$ is defined on $\ST\rest\last(p)+1$.)
  \end{defn}
 
\nin{\bf Remark} If our only aim is to obtain a model of CH \& SH, it is
enough to assume that $\gG(\gga)$ is a singleton. The assumption
that $\gG(\gga)$ is a countable collection of requirements will be used in
order to show that this forcing preserves certain Suslin trees.

  \nin{\bf Remark} If $p=(f,\gG)\in\dS$, $\gga$ is the \height\ of $p$, and
  $g:\ST\rest\gga+1\raro \Q$ is a specializing function satisfying: $\A
  x\in\dom(f)\ (g(x)\leq f(x))$, then $g$ fulfills the promise $\gG$ which $f$
  fulfills.
 
  This simple remark is used in the following way. Suppose that $p_1$ extends
  $p_0$; put $\mu_i=\last(p_i)$, $f_i=f(p_i)$ for $i=0,1$.
 
  Let $\gd$ be an order-preserving map of the set of positive rationals $\Q^+$
  into $\Q^+$ such that $\gd(r)\leq r$ for all $r$. Then define, for any
  $x\in\ST_\ga$, where $\mu_0<\ga\leq \mu_1$,
  $$
  g(x)=f_0(x\lceil\mu_0)+\gd(f_1(x)-f_0(x\lceil\mu_0)).
  $$
In words: $g$ uses $\gd$ to compress $f_1$ on $\ST\rest\mu_1+1 \setminus
\ST\rest\mu_0+1$.
 
 \nin Extend further $g$ and, for $x\in\ST\rest\mu_0+1$, define $g(x)=f_0(x)$.
  Then $(g,\gG(p_1))$ is also an extension of $p_0$ of height\ $\mu_1$.
 
  Our next aim is to show that it is possible to extend conditions to any
  height, and to enlarge promises. Then we will show properness of $\dS$. In the
  following subsection, $\dS$ is shown to specialize only those trees it must
  specialize. Then, in the next subsection $\dS$ is proved to satisfy the
  condition which allows to conclude that a countable support iteration of
  such posets adds no new countable sets.
 
  \begin{lemma}[The extension lemma]  If $p\in\dS$ and $\last(p)<\mu<\gw_1$,
  then there is an extension $q$ of $p$ in $\dS$ with $\mu=\last(q)$, and such
  that $\gG(q)=\gG(p)$. Moreover, if $h:\ST_\mu\raro\Q$ is finite and
  bounds $p$, then $h$ bounds
  an extension $q$ of $p$ of \height\ $\mu$.
  \end{lemma}
 
  \proof The `moreover' clause of the Lemma is, in fact, a direct consequence
  of the first part and the Remark above. Indeed, if $h$ bounds $p$
   as in the Lemma, pick first {\em any} extension $p_1$
  of $p$ with $\mu=\last(p_1)$, and then correct $p_1$
 as follows to obtain $q$.
 
  Put $\mu_0=\last(p),\ f=f(p)$. For some $d>0$, $\A x\in\dom(h)$
  $$
  h(x)>f(x\lceil\mu_0)+d.
  $$
 
  Let $\gd$ be an order-preserving map of $\Q^+$ into the interval $(0,d)$
  such that $\gd(x)<x$ for all $x$. Now use the Remark to correct $p_1$ and to
  obtain an extension $q$ of $p$ which satisfy for every
  $x\in\ST\rest((\mu+1)-\mu_0)$, $q(x)-p(x\lceil\mu_0)<d$. Hence $h$ bounds
  $q$.
 
  The proof of the first part of the Extension Lemma is done by induction
  on $\mu$. Since the proof is quite easy, only the outline is given.

  \nin{\bf CASE I} $\mu=\mu_0+1$ is a successor ordinal. By the inductive
  assumption, $\last(p)=\mu_0$ can be assumed, and we have to extend $f=f(p)$
  on $\ST_{\mu_0+1}$, fulfilling all the requirements in $\gG(\mu)\
  (\gG=\gG(p)$). Given any requirement $H\in\gG(\mu)$, we know that
  $H\lceil\mu_0=H_0\in\gG(\mu_0)$ is fulfilled by $f$. So, $H_0$ contains an
  infinite pairwise disjoint set of functions $h$ which bound $f$. This
  allows plenty of time to extend $f$, in $\gw$ steps, and to keep the promise
  $\gG$ at the level $\mu$.

  \nin{\bf CASE II}  $\mu$ is a limit ordinal. Pick an increasing sequence of
  ordinals $\mu_i, i<\gw$, cofinal in $\mu$. We are going to define an
  increasing sequence $p_i\in \dS$ (beginning with $p_0=p$)
 and finite
  $h_i:\ST_\mu\raro\Q$ 
  which bound $p_i$,
 by induction on $i<\gw$. Then we will set $q=(f,\Gamma)$, by
  $f=\bigcup\setm{f(p_i)}{i<\gw}\cup\bigcup\setm{h_i}{i<\gw}$, and
  $\gG=\gG(p)$.  $\last(p_i)=\mu_i$, and the passage from $p_i$ to $p_{i+1}$
  uses the inductive assumption for $\mu_{i+1}$. The role of the $h_i$'s is
  not only to ensure that $f$ is bounded on the $\mu$ branches determined by
  $\ST_\mu$, but also to ensure that the promise made in $\gG(p)=\gG$, namely
  $\gG(\mu)
$, is kept by $h=\bigcup_{i<\gw}h_i$. Each requirement
  $H\in\gG(\mu)$ must appear infinitely often in a list of missions, and at
  each step, $i<\gw$, of the definition, $h_{i+1}$ takes care of one more $h\in
  H$, so that finally an infinite pairwise disjoint subset of $H$ consists of
  functions which bound $f$.  It is here that we use the assumption that
  $\gG(\mu_i)=\gG(\mu)\lceil\mu_i$, i.e., that
  $\gG(\mu_i)=\setm{(H\lceil\mu_i)}{H\in\gG(\mu)}$.
Next we show that promised can be added.
 
  Let $p=(f,\gG)\in\dS$ be a condition of height\ $\mu$, and let $\Psi$ be
  any promise.
  We say that $p$ `includes' $\Psi$ iff for all $\gga$ such that
  $\mu\leq\gga<\gw_1$
  $$
  \Psi(\gga)\subseteq\gG(\gga).
  $$
  That is, any requirement $H\in\Psi(\gga)$ is already in $\gG(\gga)$. If $p$
  includes $\Psi$ then, obviously $p$ fulfills $\Psi$. Otherwise, it is not
  always possible to extend $p$ to fulfill $\Psi$. However, if the following
  simple condition holds, then this can be done.
 
  \begin{lemma}\label{l4.4}
[Addition of promises] Let $p\in\dS$, put $\mu=\last(p)$. Let
  $\Psi$ be a promise with $\mu<\gb=\gb(\Psi)$. Suppose for some finite
  $g:\ST_\mu\raro\Q$ (called a {\em basis} for $\Psi$), $g$ bounds $f(p)$ and
  $$
  \A \gga\geq\gb,\ \A H\in\Psi(\gga),\ \A h\in H\ (h\lceil\mu=g),
  $$
  then there is an extension $p_1$ of $p$ in $\dS$ of \height\ $\gb$ which
  includes $\Psi$.
  \end{lemma}
 
  \proof This is an easy application of the Extension Lemma. Put $f=f(p)$,
  then for some rational $d>0$, $\A x\in\dom(g)\ g(x)>f(x)+d$.
 
  Now every $H\in\Psi(\gb)$ is a dispersed collection of functions $h$ with
  $h\lceil\mu=g$. Let $p_1$ be any extension of $p$ of $\height\ \gb$; set
  $f_1=f(p_1)$. The desired extension of $p$ will be 
obtained by correcting $f_1$
  so as to fulfill $\Psi(\gb)$ and then to add $\Psi$.
This is done as follows.
 
  Let $\gd$ be an order-preserving map of the positive rationals into the
  rational interval $(0,d)$, such that $\gd(r)<r$ for every $r$.
 Define now for $x\in
  T_\ga,\ \mu<\ga\leq\gb$:
  $f_2(x)=f(x\lceil\mu)+\gd(f_1(x)-f(x\lceil\mu))$.  Then $f_2 \cup f$
 fulfills each
  $H\in\Psi(\gb)$, and thus gives the desired extension.
 
  The properness of $\dS$ is not so easy to prove, and it is here that the
  need for the promises appears. Given an elementary countable
  substructure $M\prec H(\gk)$, such that $\dS\in M$, and given a condition
  $p_0\in M$, we have to find an ``$M$-generic" condition $q$ extending $p_0$.
  In fact, we will find $q$ with a stronger property which implies that no
new reals are added: for every dense open set
  $D\subseteq\dS$ in $M$, $q\in D$. 

As in the definition of $q$ in the
  Extension Lemma (the limit case), here too an increasing sequence
  $p_i\in\dS\cap M$ of conditions and finite functions $h_i:\ST_\mu
\raro\Q$ are defined; where
  $\mu=\gw_1\cap M$. But now we are faced with an extra mission in defining
  $p_{i+1}$: to put $p_{i+1}$ in $D$, the $i$-th dense open subset of $\dS$ in
  $M$ (in some enumeration of the countable $M$). The problem with this
  mission is that perhaps whenever $r$ extends $p_i$ is in $D$, then $h_i$ does
not bound $r$.
 
  To show that this bad event never happens, requires the following main Lemma.
 
  \begin{lemma}  Let $\ST$ be an \Aron\ tree. Let $M\prec H(\gk)$ be a 
countable elementary substructure,
  where 
$\gk$ is some big enough cardinal; $\ST,\dS=\dS(\ST)\in M$. Let $p\in M$ be a
  condition in $\dS$, $\mu=\gw_1\cap M$ and $h:\ST_\mu\raro\Q$ be a finite
  function which bounds $p$. Let $D\subseteq\dS$, $D\in M$ be dense open.
  Then there is an extension of $p$, $r\in D\cap M$, such that $h$ bounds
   $r$.
  \end{lemma}
 
  \proof
  Assume for the sake of a contradiction that this is not so, and let $\ST, M,
  p,h$ etc. be a counterexample. Let $\mu_0=\last\ {p}$; $\bar{x}=\domain\ (h)$
  enumerated in some way; so  $\bar{x}\in\ ^n\ST_\mu$,
$\bar{x}=\seqn{x_0}{x_{n-1}}$. 
 Put $\bar{q}=h(\bar{x})$;
  that is, $q_i=h(x_i)$. Denote
  $\bar{v}=\bar{x}\lceil
\mu_0$; then $\bar{v}\in\ ^n\ST_{\mu_0}$, and we may assume
  $v_i\neq v_j$ for $i\neq j$ (or else, extend $p$ above the splittings of
  $\bar{x}$). In $M$: 
\medskip
 
{\it  If $r\in D\cap M$ extends $p$, then $h$ does not bound $r$.}
 
   \nin Put $g_0=h\lceil\mu_0$. Then $g_0\in M$. Say that a finite function
  $g:\ST_\gga\raro \Q$ is {\em bad} iff
\ben
\item  $\mu_0\leq\gga<\gw_1$, and
  $g\lceil\mu_0=g_0$.  
\item Whenever $r\in D$ extends $p$ and $\gga\geq\last(r)$,
  $g$ does not bound $r$.
 \een
  In other words, $g$ is bad if it mimics $h\lceil
\gga$, but it may live on other
  $n$-tuples of $\ST$. Of course, $h\lceil
\gga$ itself is bad for any $\gga$ with
  $\mu_0\leq \gga<\mu$. It follows that, in $M$ and hence in $H(\gk)$, there
  are uncountably many bad $g$'s.  Indeed, if there were only countably many
  bad functions, there would be a bound $\gga$, in $M$, for
  $\setm{\height(g)}{g\mbox{ is bad}}$; and as $\gga<\mu$, $h\lceil
\gga$ would not
  be bad.
 
  Observe that if $g$ is bad and $\mu_0\leq\gga_0<\height(g)$, then
  $g\lceil\gga_0$ is bad too.
 
  Now put
  $$
  B=\setm{dom(g)}{g\mbox{ is bad}}.
  $$
 
  Then $B$ is uncountable and closed downwards (above $\mu_0$) subset of
  $\bigcup_{\mu_0\leq\gga<\gw_1}\ ^n\ST_\gga$.  As $\ST$ is 
an \Aron\ tree, Lemma ~\ref{l1.1}
  implies that for some $\gb>\mu_0$ and some $B^0\subseteq\SB$, if we put
  $B^0_\gga=B^0\cap\ ^n\ST_\gga$, then
  \ben
  \item For $\gb\leq\gga_0<\gga_1<\gw_1$,\
  $B^0_{\gga_0}=B^0_{\gga_1}\lceil\gga_0$, and
  \item $B^0_\gb$ (and thus every $B^0_\gga$, $\gb<\gga$) is dispersed.
  \een
 
  We may find $B^0$ in $M$, since only parameters in $M$ were mentioned in its
  definition. For $\gb \leq \gga < \omega_1$,  let $\Psi(\gga)$ consists
 of $ H_\gga= \setm{g}{g\mbox{ is bad and }
dom(g)\in B^0_\gga}$.
  By Lemma ~\ref{l4.4} (Addition of promises),
 there is an extension $p_o$ of $p$ in $M$ of height $\gb$ which
  includes $\Psi$. That is, if $\gG_0=\gG(p_0)$, then for every
  $\gga\geq\last(p_0)=\gb$, 
  $H_\gga\in\gG_0(\gga)$.
 
  Now let $r\in\dS$ be {\em any} condition extending $p_0$ and in $D$.  Let
  $\gga=\last(r)$. Since $r$ fulfills $\gG_0$, for some $g\in H_\gga$, $g$ 
  bounds $r$. But this contradicts the fact that $g$ is bad.
 
  \subsection{Specialization, while Safeguarding Suslin trees}
 
  Suppose that 
we care  about a Suslin tree \SS, and wish to specialize an \Aron\
 tree
  \ST\ while keeping \SS\ Suslin. Obviously, this is not always possible: for
  example if \SS\ {\it is} \ST, or if they contain isomorphic uncountable 
subtrees. We
  will show in this section that, if \ST\ remains \Aron\ even after the addition
  of a cofinal branch to \SS, then the poset $\dS(\ST)$ specializes \ST\ while
  keeping \SS\ Suslin. 
 
  \begin{theorem} \label{t4.6}
 Let \SS\ be a Suslin tree, and \ST\ be an \Aron\ tree such that
  $\|\ST\mbox{ is \Aron}\|^{\SS}=\bone$. Then $\|\SS\mbox{ is Suslin
  }\|^{\dS(\ST)}=\bone$.
  \end{theorem}
 
  \proof  The forcing poset $\dS=\dS(\ST)$ was defined in the previous
  subsection and shown there to be proper. To prove the theorem we let $\bD$ be
  a name in $\dS$ forcing of a dense open subset of the tree $\SS$. We will
  find a condition $p\in\dS$ (extending an arbitrarily given condition in
  $\dS$) such that for some $\mu<\gw_1$, $p\force_{\dS}\ 
\SS_\mu\subseteq \bD$. This
  is enough to show that $\dS$ does not destroy the Susliness of $\SS$. The
  framework for the construction of $p$ is similar to the one for showing the
  properness of $\dS$,
  and the
  following Lemma suffices for the proof of the Theorem.
 
  \begin{lemma}  Let $\SS$ and $\ST$  be as in the Theorem.
Let $\bD$ be a name in
  $\dS=\dS(T)$ forcing of a dense open subset of  $\SS$.
Let $M\prec H(\gk)$ be a countable elementary substructure, containing
  $\ST,\SS, \bD $,
  and let $p_0\in \dS\cap M$ be a condition. 
  Let $\mu=\gw_1\cap M$, and $h_0:\ST_\mu \raro\Q$ 
be a finite function which bounds
  $p_0$. For any $b\in\SS_\mu$ there is an extension $p\in\dS\cap M$ of
  $p_0$ such that $h_0$ bounds $p$, and $p\force_{\dS}$ $b\in \bD$.
  \end{lemma}
 
  \proof Assume that this Lemma does not
  hold. Let $M, p_0, h_0$ etc. be a counterexample. Put $\mu_0=\last(p_0)$,
  and $g_0=h_0\lceil\mu_0$.
 
  The Suslin tree $\SS$ is a c.c.c. forcing notion which adds no new countable
  sets. We are going to define first a {\em name} $\bB$, in $\SS$ forcing, of an
  uncountable tree of `bad' functions, and derive a promise $\gG$ out of this
  $\bB$, a promise which, when adjoined to $p_0$, will give the desired
  contradiction.
 
  Forcing with $\SS$, extend $b\in\SS_\mu$ to a (generic) branch $G$ of $\SS$,
  and let $V[G]$ be the extension of the universe $V$ thus obtained. We have:
  $$
  M[G]\prec H(\gk)[G]=H(\gk)^{V[G]}.
  $$
 
  In $V[G]$, and hence in $M[G]$, $\ST$ is still an \Aron\ tree by the
  assumption of the Theorem. The following definition is carried out in $V[G]$,
  but all its parameters are in $M[G]$:
 
  \begin{defn} A finite function $h:\ST_\gga\raro\Q$ is bad iff:\\
  \nin 1. $\mu_0\leq\gga<\gw_1$, and $h\lceil\mu_0=g_0$.\\
  2. Whenever $p\in\dS$ extends $p_0$ and $\gga\geq \last{p}$ and $G_\gga=e$, if
  $p\force_{\dS}\ e\in \bD$ then $h$ does not dominate $p$.
  (Recall that $G_\gga$ is the unique node in $G\cap \SS_\gga$.)
\end{defn}

 For any
  $\mu_0\leq\gga <\mu$, $h_0\lceil
\gga$ is bad. (If not, by elementarity of $M[G]$,
  there is, in $M$, an extension $p$ of $p_0$, of height $\gga$ and such that
  $h_0\lceil\gga$ bounds $p$ and $p\force_{\dS}\ (G_\gga)\in \bD$. But 
then,
  as $b>G_\gga$, $p\force_{\dS}\ b\in \bD$, in contradiction to our
assumption.) Hence the set of bad functions is uncountable.
 
  Obviously, if $h$, of \height\ $\gga$, is bad and $\mu_0\leq\gga'<\gga$,
  then $h\lceil\gga'$ is bad too.
 
  We know how to find  (in $M[G]$) an ordinal $\mu_0 \leq \gb < \go_1$,
and a collection $B(\gga)$, $\gb\leq\gga<\gw_1$,
 such that  $B(\gga)$ is a set of bad functions of 
height $\gga$, and\\
  1. For $\gb\leq\gga_0<\gga_1<\gw_1$, $B(\gga_0)=B(\gga_1)\lceil\gga_0$,\\
  2. $B(\gb)$ is dispersed.
(See Definition ~\ref{d4.1} (4)(b), and Lemma ~\ref{l1.1}.)

  Let $\bB\in M$ be a name of $B$ in $V^{\SS}$,
 and let $b_0<b$ be a condition in
  $\SS$ which forces these properties of $B$. In particular, $b_0$ forces
  ``all functions in $\bB(\gga)$ are bad".
 
  Now, back in $V$, we define the promise $\gG$. For every countable
  $\gga\geq\gb$, $\gG(\gga)$ is the collection of all requirements $H$ of
  height $\gga$ such that $\|H=\bB(\gga)\|^{\SS}>{\mbox{{\small\bf
$0$}}}$. Again $\gG\in M$. Since
  $\SS$ is a c.c.c. poset, $\gG(\gga)$ is countable, and since $\SS$ adds no
  new countable sets, $\gG(\gga)$ is non-empty (some condition in $\SS$ above
  $b_0$ `describes' $\bB(\gga)$) and $\gG(\gga)$ is countable.
 
  Since $g_0$ is a basis of $\gG$, and $g_0$ bounds $p_0$, there
  is an extension $p_1$ of $p_0$, in $\dS\cap M$, which includes $\gG$.
(See the Addition of Promises Lemma ~\ref{l4.4}.)
 
  Next, find a node $d\in \SS$, with
 $b_0<d$, such that $p_2\force_{\dS }\ d\in \bD$ for some extension $p_2$ of
$p_1$ with $p_2\in\dS\cap M$.
 This is possible since $\bD$
 is assumed to be a
  name in $V^{\dS}$ such that $\|\bD\mbox{ is dense in } \SS\|^{\dS}=\bone$.
 
  Let $\gga=\last (p_2)$, and let $d_1>d$ be a node in $\SS$ which forces
  ``$H=\bB(\gga)$" for some requirement $H$ of height $\gga$. Then
  $H\in\gG(\gga)$ and so some $h\in H$ bounds $p_2$ (as $p_2$
  fulfills $\gG$).  But $d_1\force_{\SS}\ h$ 
{\em is bad}, contradicts $p_2\force_{\dS}\ 
  d\in \bD$.
 
  \section{$\ga$-properness and $\gw_1$-\hbox{D \kern -20pt D}-completeness 
of \dS}\label{s5}
 
In chapter V of Shelah [1982] (Sections 3,5 and 6) the notions of
$\ga$-properness and $\dD$-completeness are defined, and an $<\go_1$-proper,
simple $\dD$-complete forcing which specializes an \Aron\ tree is described.
Section 7 there shows that the iteration of such forcings adds no reals.
In chapters VII and VIII different notions of chain conditions
are introduced: the $\aleph_2$-e.c.c and the  $\aleph_2$-p.i.c. Any of them
can be used to show that our iterations satisfy the $\aleph_2$-c.c. (The
second is particularly useful if $2^{\aleph_1}>\aleph_2$). We shall review
here these definitions, but will not give proofs for the preservation
theorems which may be found in the Proper Forcing book.

  To use the theory of proper forcings which add no reals we will show that
 (1) the specializing poset $\dS=\dS(\ST)$ is $\ga$-proper for every
  $\ga<\gw_1$, and that (2) for some simple $\gw_1$-completeness system \dD,
  \dS\ is \dD-complete.
 
  By ``a tower of \length\ $\ga+1$ of substructures of $H(\gl)$" we mean here a
  sequence $\bar{N}=\sseqn{N_i}{i\leq \ga}$ of countable $N_i\prec H(\gl)$ such
  that
  \ben
  \item $\bar{N}$ is continuously increasing. ($N_\gd=\bigcup_{i<\gd} N_i$,
for limit $\gd \leq \ga$).
  \item $\sseqn{N_j}{j\leq i}\in N_{i+1}$.
  \een

  \begin{defn} $P$ is $\ga$-proper ($\ga<\gw_1$), iff for every large enough
  $\gl$ and tower $\sseqn{N_i}{i\leq\ga}$ of countable 
substructures of $H(\gl)$ of
  \length\ $\ga+1$, if $P\in N_0$ and $p\in P\cap N_0$, then there is an
  extension $q$ of $p$ in $P$ such that $q$ is an $(N_i,P)$-generic condition
  {\em for every} $i\leq\ga$.
  \end{defn}
 
  \begin{theorem} 
 $\dS$ is $\ga$-proper for every $\ga<\gw_1$.
  \end{theorem}
 
  We only {\em indicate}
 the proof since there is not much to say which was not said
  for the case $\ga=1$. The proof is by induction on $\ga$. The case of a
  successor ordinal is an obvious application of the inductive assumption and
  the properness of $\dS$. In case $\ga$ is a limit ordinal, given a tower
  $\sseqn{N_i}{i\leq \ga}$ as in the definition, pick an increasing
  $\gw$-sequence $i_n<\ga$, $n<\gw$, converging to $\ga$. We will define an
  increasing sequence $p_n\in N_{i_{n}+1}$ such that $p_n$ is
  $(N_j,P)$-generic condition for every $j\leq i_n$. The inductive assumption
  and the assumption that $\sseqn{N_k}{k\leq i_n}\in N_{i_{n}+1}$ are used to
  get $p_n$ in the elementary substructure $N_{i_{n}+1}$. We must be careful
  so that $f\stackrel{\rm def}{=}
\bigcup_{n<\gw} f(p_n)$ is bounded on every branch of $\ST\rest\ga$
  determined by points in $\ST_\ga$, and that $f$ fulfills every requirement in
  $\gG(\ga)$ for $\gG=\gG(p_n),\ n<\gw$.
  But we learned how to do it when proving properness of \dS.
 
  The \dD-completeness of \dS\ is equally simple, if only the definition is
  clear. Let us review it on the informal level first.
 
  Think on the difference between the poset \dP, for adding a new subset to
  $\gw_1$ with countable conditions on the one hand, and the poset $\ST$, a
  Suslin tree, on the other hand. Both posets add no new countable sets, but
  while posets like \dP\ can be iterated without adding reals, an
  iteration of Suslin trees can add a new real (see Jensen and 
Johansbr\aa ten [1974]).
 
  Pick a countable $M\prec H(\gl)$ with $\dP, \ST\in M$ and look for \dP-generic
  and $\ST$-generic filters $G_{\dP}$ and $G_\ST$ over $M$ which have an upper
  bound in \dP, and in $\ST$ (this is what it takes to show ``no new countable
  sets are added").  While 
$G_{\dP}$ can be defined, in a sense, from within $M$;
  the definition of $G_
\ST$ requires knowing $\ST_\ga$.  To clearly see this
  difference, let $\Pi:M\raro\ovl{M}$ be the transitive collapse of $M$ onto
  the transitive structure $\ovl{M}$. If we only have $\ovl{M}$ at hand (and a
  countable enumeration of $\ovl{M}$) 
then we can define a \dP-generic filter over
  $\ovl{M}$, and any such filter has an upper bound in \dP. However, for 
$\ST$
  the situation is radically different: even though any branch of $\ST\cap M$ is
  $\ovl{M}$-generic, there is no way to know which branches have an upper
  bound in $\ST$, unless $\ST_\ga$ is given to us.
 
  For the poset $\dS(\ST)$ ($\ST$ now an \Aron\ tree) the situation is subtly in
  between \dP\ and $\ST$: It seems that 
we need to know $\ST_\ga$ (and more) to define $M$-generic filters
  over \dS, but in fact 
this is less crucial: there is room for some errors. Let us
  make this more precise in the following. Recall the properness
  proof, and suppose that $M\prec H(\gl)$  and the collapse $\Pi:\ M \raro
\ovl{M}$
are given. We seek to find a generic filter $G$ over $\ovl{M}$ such that
$\Pi^{-1}G$ has an upper bound in $\dS$.
   Besides $\ovl{M}$, the only
  parameters of importance were $\ST_\mu$
  ($\mu=\gw_1\cap M$) and the function $\gG$ which assigns to each $p\in\dS\cap
  M$, the countable set $\gG(p)(\mu)$ of requirements (in the sense of
$\ST'_\mu$) of height $\mu$.
Suppose that not the real $\ST_\mu$ and $\gG$ are given, but just a
countable set of cofinal branches of $\ST \cap \ovl{M}$ (called $\ST'_\mu$) and
any function $\gG'$ such that $\gG'(p)$ is a countable collection of
requirements of height $\mu$, and
for every $p\in\dS\cap M$ and
  $\gb\in\gw_1\cap M$, $\gG'(p)(\gb)=\setm{X\lceil \gb}{X\in \gG'(p)}$.
Then the increasing, $\ovl{M}$-generic sequence of conditions,
$p_i\in\dS\cap \ovl{M}$, 
  $i<\gw$, could be defined to give a filter $G$.  Of
  course, if $\ST'_\mu$ and $\gG'$ are arbitrary, then we cannot be sure that 
$\Pi^{-1}$ of the
  filter $G$ thus obtained has an upper bound in \dS. {\em However},
  the following observation comes to our rescue: given a countable collection
  $\setm{\pair{\ST^i_\mu}{\gG^i}}{i<\gw}$, it is possible to find a generic $G$
  which is good for {\em every}
 $\pair{\ST_\mu^i}{\gG^i}$. `Good' in the sense that if some $\pair{\ST^i_\mu}
{\gG^i}$ were the real thing,
  then $G$ would have an upper bound in the external \dS. This is the essence
  of the notion of simple \dD-completeness. For completeness, we give now
the definition from chapter V of
  Shelah's Proper Forcing [1982]. The reader can then 
see that \dS\ is indeed \dD-complete for an
  $\omega$-completeness simple system.
The theory developed there shows that the countable support iteration of
$\ga$-proper ($\ga < \gw_1$) and simple
\dD-complete posets adds no new countable sets.

\nin{\bf Definition of \dD-completeness}. (See Definitions 5.2, 5.3 and 5.5
in Shelah [1982]). For any structure $N$, let $\pi: N \rightarrow \bar{N}$
denote the Mostowski collapse of $N$ to a transitive structure. When enough
set-theory is present in $\bar{N}$, the forcing relation can be defined in
$\bar{N}$. So, given a poset $\barP \in \bar{N}$, if $\bar{N}$ is countable,
an $\barN$-generic filter over $\barP$ can be found, and the generic
extension $\barN[G]$  can be formed. We let\\
$Gen(\barN,\barP,\barp)=\{G\subset \barP \mid G \ \mbox{\it is } \barN-{\it
generic}\ \mbox{\it filter over}\ \barP,\ {\it and}\ \barp \in G\}.$

The function $\dD$ is called an $\aleph_1$-completeness system iff\\
{\it For every countable transitive model} $\barN$  {\it (of enough set-theory)
and} $\barp \in \barP \in \barN$, $\dD(\barN, \barP, \barp)$  
{\it is a family of subsets of} $Gen(\barN, \barP, \barp)$
 {\it such that every
intersection of countably many sets in that family is non-empty.}

Thus if $G \in A \in \dD(\barN,\barP,\barp)$ then $G$ is an $\barN$-generic
filter over $\barP$ with $\barp \in G$, and if $A^i \in 
\dD(\barN,\barP,\barp)$ then $\cap_{i<\gw}A^i$ is non-empty.

Given a completeness system $\dD$, we say that the poset $P$ is
$\dD$-complete iff for some large enough $\gk$ the following holds: For
every $N \prec H(\gk)$, with $\barP \in N$ and for every $p \in P$, let $\pi:N
\rightarrow \barN$ be the transitive collapse of $N$; put $\barP=\pi(P)$,
$\barp=\pi(p)$. There is some $A \in \dD(\barN,\barP,\barp)$
such that for every $G \in A$:\\
$\pi^{-1}(G)=\{\pi^{-1}(g) \mid g \in G \}$ contains an upper bound in $P$.

Finally, let us say that the completeness system $\dD$ is {\em simple} iff
it is given by a formula $\psi(G,\barP,\barp,x)$ in the following way:\\
$\dD(\barN,\barP,\barp)=\{A_x \mid x \subset \barN\}$, where\\
$A_x=\{G \in Gen(\barN,\barP,\barp)\mid \langle \barN \cup \CP(\barN),\in
\rangle \models \psi(G,\barP,\barp,x) \}.$ 

In our case, the parameter $x$ describes $\ST_\ga$ and the function $p
\longmapsto \gG(p)(\ga)$, where $\ga=\gw_1^{\barN}$.

As for the $\aleph_2$-chain condition of $\dS$, it follows from CH by the
obvious remark that if two conditions have the same specializing function
(but different promises) then they are compatible. In Chapter 8 of Shelah
[1982] the notion of $\aleph_2$-p.i.c ($\aleph_2$ proper isomorphism
condition) is defined, and it is shown that countable support iteration
of length $\gw_2$
of such posets satisfies the $\aleph_2$-c.c. if CH is assumed. Our posets
$\dS$ clearly satisfy the $\aleph_2$-p.i.c. and hence that result may be
applied to conclude that the $\aleph_2$ chain condition holds for the
iteration.
  \section{Models with few Suslin trees} \label{s6}
 
  Suppose $\SS$ and all the  derived trees of $\SS$ are
  Suslin trees. We shall
 find now  a generic extension in which the only Suslin trees are $\SS$ and
  its derived trees, and such that no new countable sets are added by this
  extension. The extension is obtained as an iteration of length
$2^{\aleph_1}=\aleph_2$ of
  posets of type $\dS(\ST)$ described as follows.
 
  By the result of Section ~\ref{s3},
 we know that if $\SS$ and its derived trees
  remain Suslin at each stage of the iteration,
 then this also holds for 
the
 final limit of the iteration. We know that no new countable sets are added
by the iteration of $\dS(\ST)$ forcings,
  and that the $\aleph_2$-chain condition holds. The definition of the
  iteration is such that if,  for some $\ga<\gw_2$, $\dP_\ga$ is defined, then
  $\dP_{\ga+1}$ is obtained as an iteration $\dP_{\ga}*\dqq_\ga$ where (in
  $V^{\dP_\ga}$) $\dqq_\ga$ has the form $\dS(
\ST)$ for an `appropriate' \Aron\ tree ($\ST$ is appropriate if
  for any derived tree $\SS^1$ of $\SS$, $\dline{\ST\mbox{ is
  \Aron}}^{\SS^1}=\bone$). Then, by Theorem ~\ref{t4.6},
 $\dline{\mbox{all derived trees
  of } \SS \mbox{ are Suslin }}^{\dS(T)}=\bone$.
 
  When care is taken of all appropriate
 $\ST$ as above, the final extension $V[G]$ satisfies
  \ben
  \item $\SS$ and its derived trees are Suslin.
  \item For any \Aron\ tree $\ST$, either $\ST$ is special, or  for some derived
  tree $\SS^1$ of $\SS$, $\dline{\ST\mbox{ is not \Aron}}^{\SS^1}>{\bf 0}$.
  \een
 
  The latter possibility implies, as the following Lemma shows, that $\ST$
  contains a club-isomorphic copy of a derived tree of $\SS$.
 
  \begin{lemma} Assume that
 $\SS$ and its derived trees are all Suslin. Suppose $\ST$
 is an \Aron\ tree, and $\dline{\ST\mbox{ is not \Aron}}^{\SS^{1}}=\bone$ for a
  derived tree $\SS^1$ of $\SS$. Assume further that $\SS^1$ is of least 
dimension
  with this property. Then $\SS^1$ is embeddable on a club set into $\ST$.
  \end{lemma}
 
  \proof 
  Let us fix first some notation. $\SS^1$ has the form
  $\SS_{a_{1}}\times\ldots\times \SS_{a_{n}}$ for some $n$-tuple $\seqn{a_1}{a_n}$
of distinct elements of $S_\gga$ for some $\gga<\gw_1$. $n$ is called the
dimension of $\SS^1$.
  Now, when we say that $e\in \SS^1$ is of the form $e=\seqn{e_1}{e_n}$ it is
  assumed that $e_i>a_i$ in $\SS$.

For this Lemma, we assume that any node of limit hight in $\ST$ is
determined by its predeccessors.
Let $b$ be a name in $V^{\SS^{1}}$ such that $$\dline{b\mbox{ is a
  cofinal branch in }\ST}^{\SS^{1}}=\bone.$$ 
For any $e_1\in \SS^1$ and $\ga<\gw_1$
  there is an extension $e_2$ of $e_1$ which determines $b_\ga=b\cap \ST_\ga$.
  That is, for some $x\in \ST_\ga$, $e_2\force_{\SS^1}\
 x\in b$. The set $D_\ga$ of
  all conditions in $\SS^1$ which thus determine $b_\ga$ is a dense open set in
  $\SS^1$. Since $\SS^1$ is Suslin, there is a club set $E\subseteq \gw_1$ of
  limit ordinals such that for $\ga\in E$, if $e\in (\SS^1)_\ga$ then $e\in
  D_\gb$ for all $\gb<\ga$. That is, $e$ determines $b_\gb$ for all $\gb<\ga$.
  But then $e$ must determine $b_\ga$ as well, since there is a {\em single}
  node in $\ST_\ga$ above all those determined $b_\ga $'s. So, we have that for
  $\ga\in E$, $(\SS^1)_\ga\subseteq D_\ga$; hence for every $e\in(\SS^1)_\ga$ there is
  some $f(e)\in \ST_\ga$ with $e\force_{\SS^1}\
 b_\ga=f(e)$. As we will see, $f$
  is an embedding of $\SS^1$ on some club into $\ST$. Clearly $f$ is an order
  preserving map of $\SS^1\rest 
E$ into $\ST\rest E$.  We will find a club set $D\subseteq
  E$ such that, on $\SS^1\rest 
D$, $f$ is one-to-one. The basic observation is that,
  as $\ST$ is \Aron\ (and any node in $\SS^1$ has extensions to every higher
  level), every $e\in \SS^1\rest E$ has
 two extensions, $e_1$ and $e_2$ such that
  $f(e_1)\neq f(e_2)$. The following is a slight strengthening of this,
which is obtained from
  the minimality of the dimension $n$ of $\SS^1$.
 
  \nin{\bf Claim:} For every $e\in \SS^1$, and for every set of indices
  $h\subset\fsetn{1}{n}$ (strict inclusion) there are two extensions, $e'=\seqn{e'_1}{e'_n}$ and
  $e''=\seqn{e''_1}{e''_n}$ of $e$, such that $f(e')\neq f(e''
)$; and
  $e'_i= e''_i$ for $i\in h$.
 
  \proof Suppose this is not so, and for some
  $h=\fsetn{h(1)}{h(k)}\subseteq\fsetn{1}{n}$ with $k<n$, for some $e=\seqn
  {e_1}{e_n}$ for every two extensions $e'$ and $e''$ of $e$ with the same
  restriction on $h$, $f(e')=f(e'')$. This means that restricted to $\SS^1_e$,
  the function $f$ actually depends on $\SS^2=
\SS_{e_{h(1)}}\times\ldots\times
  \SS_{e_{h(k)}}$. This enables us to define a name, in $\SS^2$ forcing, of a branch
  in $\ST$. But the dimension $k$ of $\SS^2$ contradicts the minimality of $n$.

Now the proof of the Lemma can be concluded by showing that the embedding
$f$ defined above is one-to-one on a club set. This follows from the Claim
since not only $\SS^1$ but any other derived tree of $\SS$ is Suslin. Take
for example a countable elementary substructure, $M$, of some $H_\gk$, and
put $\delta=M \cap \go_1$. We claim that if $e^1 \not = e^2$ are in
$(\SS^1)_\gd$ then $f(e^1 \not = f(e^2)$. Confusing sequences with sets,
put $e=e^1 \cup e^2$; then for some $k$ with 
$n<k \leq 2n$, $e$ is a $k$-tuple. $e$
`is' in fact an $M$-generic branch of a derived tree of $\SS$ of dimension
$k$. What the Claim implies is that for a dense open set of $k$-tuples of the
form $e' \cup e^{''}$ in that derived tree, $f(e') \not = f(e^{''})$. Since
$e$ is in that  dense set,
   $f(e') \not =f(e^{''})$. 

The generalization of our discussion
 to any collection of trees poses no problems. Suppose we are given a
  collection $\CU$ of Suslin trees such that if $\SS\in \CU$ then all derived
  trees of $\SS$ are Suslin. Assume: $2^{\aleph_0}=\aleph_1$ 
and $2^{\aleph_1}=\aleph_2$. 
   Then iterate $\dS(\ST)$ posets, just as before, so that all trees in
  $\CU$ and their derived trees remain Suslin. We know that  this is
  possible, for any \Aron\ tree $\ST$, unless $\dline{\ST\mbox{ is not
  \Aron}}^{\SS^{1}}>0$ for some $\SS^1$ which is a derived tree of a tree in 
$\CU$.
  In such a case, we know that $\ST$ must contain a restriction to a club set of
  a derived tree of some $\SS\in \ST$.
  \section{The uniqueness of simple primal Suslin sequences}\label{s7}
 
  This section sets the preliminaries needed to prove the main theorem: the
  notions of {\em simple} and {\em primal}
 sequences of Suslin trees are defined, and
  the uniqueness of such sequences is proved. Using this material and the 
machinery
  developed to construct Suslin trees and to specialize them at will, the
  Encoding Theorem will be easily demonstrated in the following section.

We will deal here not only with $\gw_1$-sequences of \Aron\ trees, but also
with $I$-sequences, $\CT=\sseqn{\ST^\zeta}{\zeta\in I}$, of \Aron\ 
trees, where $I$ is an $\gw_1$-like set of indices.
 
  A linear order $(I,<)$ is said to be $\gw_1$-like iff it is uncountable but
  all proper initial segments are countable. In this paper we need a slightly
  stronger version, and add to these requirements of $\gw_1$-like that any
  point has a successor, and that a first element exists.
 
  We say that $a\in I$ is a `limit' point if it is not a successor (so the
  first element is a limit). A point $a\in I$ is said to be `even' iff it is a
  limit or it has the form $\gd+n$, where $\gd$ is a limit and $n<\gw$ is an
  even integer. Similarly
   `odd' points of $I$ are defined.
 
  We will call the members of $I$ `indices', since this is how they will be
  used. In some cases, $I$ is  or is isomorphic to
  $\gw_1$, but in general an $\gw_1$-like order need not to be well-founded.
  Indeed, an important point of our argument is that, in some universe, the
  Magidor-Malitz quantifiers can force $I$ to be well-founded.

  Let $\CT=\sseqn{\ST^\zeta}{\zeta\in I}$ be a sequence of \Aron\ trees.
  Recall
  that for $d\in[I]^{<\gw}-\{\emptyset\}$ (a finite non empty subset of
  $I$), $\ST^d=\bigcup_{\zeta\in d}\ST^\zeta$ is the disjoint union of the
  \Aron\ trees with indices in $d$. A derived tree of $\ST^d$ is thus a
  product of derived trees of the $\ST^\zeta$'s. One of these derived trees is
  $\ST^(d)=\bigtimes_{\zeta\in d}\ST^\zeta$.
 
  Let \st\ be a collection of non-empty finite subsets of $I$ which is closed
under subsets, and let
  $\spa=[I]^{<\gw}-\{\empt\}-\st$ be the complement of \st.
  (\spa\ is closed under supersets.) We say then that
  $(\st,\spa)$ is a pattern (over $I$)

  \begin{defn} We say that the $I$-sequence $\CT$ of \Aron\ 
trees has the pattern
  $(\st,\spa)$ if
  \ben
  \item For $d\in\st$, every derived tree of $\ST^d$ is Suslin.
  \item For $d\in\spa$, $\ST^{(d)}$
 is a special tree.
  \een
  \end{defn}
 
  \begin{defn} \ben
  \item A collection $\CU$ of Suslin trees is {\em primal} iff all derived trees of
  trees in $\CU$ are Suslin, and for any Suslin tree $\SA$ there exists some
  $\SS\in\CU$ such that a derived tree of $\SS$ is club-embeddable into $\SA$.
  \item The sequence $\CT$ with Suslin-special pattern $(\st,\spa)$ is
  called primal iff the collection $\CU=\setm{\ST^d}{d\in\st}$ is primal.
  \een
  \end{defn}
 
  We may summarize our results obtained so far in the following:

  \begin{theorem} 
 {\rm 1.} Assume $\diamond_{\gw_{1}}$. Given any pattern 
$(\st,\spa)$ over $\gw_1$,
 there exists an $\gw_1$-sequence of \Aron\ trees
  $\CT=\sseqn{\ST^\zeta}{\zeta\in\gw_1}$ with this pattern. (Section
~\ref{s2.3}.)\\
  \nin{\rm 2.} Assume $2^{\aleph_0}=\aleph_1$ and $2^{\aleph_1}=\aleph_2$.
Let $\CU$ be
  a collection of Suslin trees such that for $\SS\in\CU$ all derived trees of
  $\SS$ are Suslin as well.
  There is then an $\aleph_2$-c.c. generic extension which adds no new
  countable sets, and in which $\CU$ is a primal collection of Suslin trees.
(Section ~\ref{s6}.)
\end{theorem}
 
  \subsection{Simple Patterns}\label{s7.1}

 Let $I$ be an $\gw_1$-like order.
  We will have to refer to quadruples $\zeta_1<\zeta_2<\zeta_3<\zeta_4$ of
  indices given in increasing order in $I$, with some simple properties
  called types.  For example, $\bar{\zeta}=\seqn{\zeta_1}{\zeta_4}$ is of type
  $\langle$odd, even, odd, even$\rangle$ if $\zeta_1$ and $\zeta_3$ are odd, and
  $\zeta_2,\ \zeta_4$ are even.  Similar notations are obvious. 
 
  Now we define when the pattern pair $(\st,\spa)$ is said to be {\em simple}:
  If \st\ contains all tuples in the columns Suslin, and \spa\ 
contains all tuples
  in the column Special, in the following table.
  In case $I=\gw_1$,
for each limit ordinal $\gd$, pick a canonical well-order $E_\gd$ of
  $\gw$ of order-type $\gd$. By standard encoding, we assume that
  $E_\gd\subseteq \gw$.
\medskip

  \begin{tabular}{ll}
  \underline{Suslin} & \underline{Special}\\ \\
All triples, pairs, and singletons & All quintuples\\
  \lan even, even, even, even\ran & \\
  \lan odd, odd, odd, odd\ran \\
  \lan odd, even, odd, even\ran& \lan even, odd, even, odd\ran\\
  \lan odd, even, even, odd\ran & \lan odd, even, even, even\ran \\
  \lan limit $\gd$, even, odd, odd\ran & \lan non-limit even, even, odd,
  odd\ran \\
  \lan $\zeta$, $\zeta+1$, odd, odd\ran & \lan $\ga$, $\gb$,odd,
  odd\ran where $\ga$ and $\gb$ have\\
 & different parity, and $\ga +1 < \gb$.\\
  \lan $\gd+1,\ \gd+3,\ \gd+4,\ \gd+(5+i)$\ran &
  \lan $\gd+1,\ \gd+3,\ \gd+4,\ \gd+(5+i)$\ran \\
 where $i \in E_\gd$, $\gd$ a
  limit &
where $\gd$ is limit and $i\not\in E_\gd$.
  \end{tabular}
\medskip
 
  In the table we used the sets $E_\gd\subseteq\gw$; these are required only
  in case $I$ is well-ordered.
 
  It should be checked that if $d\subseteq I$ appears in the Suslin column
  of the table, and $e\subseteq I$ appears in the special column, then
  $e\not\subseteq d$. The Suslin tuples  are closed under subsets,
and even after closing the
 the special tuples under supersets, disjoint sets are obtained.
 
  The following theorem explains the use of simple sequences (sequences of
trees
  with simple patterns).  In fact, as the reader may find out, there are 
  other notions of simplicity which can be used to derive the conclusion of the
  Theorem.
 
  \begin{theorem}[Unique Pattern] 
\label{t7.4}
Suppose that
  $\CT=\sseqn{\ST^\zeta}{\zeta\in \gw_1}$, and $\CA=\setm{\SA^\zeta}{\zeta\in I}$
  are $\gw_1$ and $\gw_1$-like sequences of Suslin trees with simple Suslin-special
  patterns. IF $\CT$ is primal, then $I$ is isomorphic to $\gw_1$, and $\CT$
  and $\CA$ have the {\em same} Suslin-special pattern.
  \end{theorem}
 
  \proof We are going to define an order 
isomorphism $\fnn{d}{I}{\gw_1}$, and show
  that for every $\zeta\in I$, $\SA^\zeta$ contains a club-embedding of a
  derived tree of $\ST^{d(\zeta)}$. This suffices to derive the equality of
  the patterns of $\CT$ and $\CA$. We shall use the following easy observations:
If $\SA$ is special and \SB\ is any tree, then
  $\SA\times\SB$ is special too.
  If $\fnn{h}{\SA}{\SB}$ is a club embedding of the tree $\SA$ into $\SB$,
  then \\
  \begin{enumerate}
  \item  $\SB$ is special $\Raro$  \SA\ is special (on a club set of
  levels and hence on all levels. See Devlin and Johnsbr\aa ten [1974]).
  \item  $\SA$ is special $\Raro$  $\SB$ is not Suslin.
  \end{enumerate}

  Let $(\st_1,\spa_1),\ (\st_2,\spa_2)$ be the patterns of $\CT$ and $\CA$
  respectively.
  Since $\CT$ is assumed to be primal, for every Suslin tree $\SA^\zeta$
  there is $d=d(\xi)\in\st_1$ such that $\SA^\xi$
contains a club image of
  a derived tree of $\ST^d$. Our aim is to prove that 
  $$
   d(\xi)\mbox{\it is a singleton, and } d\ 
\mbox{\it establishes an isomorphism of }I\mbox{\it onto}\ \gw_1.
  $$
 
  This will be achieved in the following steps.
  \bd
  \itm{a} {\em There is no quintuple $\xi_1,\ldots,\xi_5$ such that, for
  all indices $i,j$, $d(\xi_i)=d(\xi_j)$.} Suppose, for the sake of a
  contradiction, that for some $d\in\st_1$, $d=d(\xi_i)$ for five indices
  $\xi_1,\ldots,\xi_5$. Then there are club embeddings, $h_i$, from
  derived trees of $\ST^d$ into $A^{\xi_i}$, $1\leq i\leq 5$.
  We may combine these embeddings into a club embedding of a derived tree of
  $\ST^d$ into $\bigtimes _{1\leq i\leq 5} \SA^{\xi_i}$. But since the
  product of a quintuple of trees is a special tree, this derived tree of 
$\ST^d$ cannot be Suslin.
  \itm{b} {\em There are not uncountably many $\xi$'s with $|d(\xi)|>1$.}
  To see this, suppose the contrary, and let an uncountable set
  $X\subseteq I$ be  such that for $\xi\in X$, $d(\xi
)$ contains more than one
  element. We may assume that the finite sets $d(\xi), \xi\in X$, form a
  $\Delta$-system, and that either all members of $X$ are odd or all are even.
  Hence, for all quadruples $d \subset X$, $\SA^d$ is Suslin.
 We will get the contradiction be
  considering the two possibilities for the $\Delta$-system. If the core of
  the system is all of $d(\xi)$, i.e., $d=d(\xi_1)=d(\xi_2)$ for
  $\xi_1,\xi_2\in X$, then a contradiction to (a) is obtained.
 
  If the core of the system is strictly included in $d(\xi),$ for $ \xi\in X$,
  then for a quadruple $\xi_1,\ldots,\xi_4$ in $X$,
  $d=\bigcup_{1\leq i\leq 4} d(\xi_i)$ contains $\geq 5$ indices. Now
  $\bigtimes_{i\in d}\ST^i$ is a special tree, and has a derived tree which is
  club embeddable into $\bigtimes_{1\leq i\leq 4}\SA^{\xi_i}$ which is a
  Suslin tree. This is clearly impossible since a derived tree of a special
  tree is special. 
 
  Now that we have proved that on a co-countable set $d(\xi)$ is a
  singleton, we proceed to show that $d(\xi)$ is a singleton for every
  $\xi$.
  \itm{c} {\em For every $\xi\in I,\ |d(\xi)|=1$.} This will enable us to
  change notation and write $d(\xi)\in \go_1$ (instead of $d(\xi)\subseteq
  \go_1$). Assume, for some $\xi_1$, $|d(\xi_1)|\geq 2$. Suppose, for example,
  that $\xi_1$ is even. We can find (in the co-countable set of (b)) even
  indices $\xi_2,\xi_3,\xi_4$ such that $c=\bigcup_{1\leq i\leq 4}
  d(\xi_i)$ contains $\geq 5$ indices. Since
  $e=\fsetn{\xi_1}{\xi_4}\in\st_2$, all derived trees of $A^e$ are Suslin,
  and in particular $B=\bigtimes_{i\in e}A^i$ is Suslin. On the other hand,
  there is an embedding of a derived tree of $\ST^c$ on a club into \SB,
  but this is impossible as any such derived tree of $\ST^c$ is special (as
$\mid c \mid \geq 5$).
  \itm{d} $d$ {\em is one-to-one}. Suppose that $\xi_1<\xi_2$ but
  $d(\xi_1)=d(\xi_2)$. Consider the four possibilities for
  $(\xi_1,\xi_2)$. (1) both are even, (2) both are odd, (3) $\xi_1$ is
  even and $\xi_2$ is odd, (4) $\xi_1$ is odd and $\xi_2$ is even.
 
  In each case, it is possible to find $\gz_1,\ \gz_2$ such that
  $e=\{\xi_1,\xi_2\}\cup \{ \gz_1,\gz_2\} \in\spa_2$. (For example, if both
$\xi_1$ and $\xi_2$ are even, find odd  $\gz_1,\ \gz_2$ such that
$\xi_1<\gz_1 < \xi_2 < \gz_2$.) Yet
  $t=\setm{d(\xi)}{\xi\in e}$ contains at most 3 indices, and a club embedding
  of a derived tree of $\ST^t$ into $\bigtimes_{\xi\in e}\SA^\xi$ is
  obtained. This contradicts the fact that the first tree is Suslin and the
  latter is special.

At this stage we don't know yet that $d$ is order preserving; but as $d$ is
one-to-one from an $\gw_1$-like order into an $\gw_1$-like order, for any
$\ga$, if $\gb > \ga$ is sufficiently large, then $d(\gb)>d(\ga)$. This
simple remark is used below.
  \itm{e} $\xi$ {\em is even if and only if $d(\xi)$ is even}.
  Let us first check this: Could it be that there are both 
 uncountably many even
  $\xi$'s with $d(\xi)$ odd, and uncountably many odd $\xi$'s with
  $d(\xi)$ even? No, because in such a case (using the remark made above)
 we will find a quadruple
  $\seqn{\xi_1}{\xi_4}$ of type \lan even, odd, even, odd\ran\ with
  $d$-image of type \lan odd, even, odd, even\ran. But this is impossible
  since the first type is in the special column, and the second in the Suslin
  column of the simplicity table.
 
  It follows now that there are not uncountably many odd $\xi$'s for which 
$d(\xi)$
  is even. For otherwise, using our result above, on a co-countable
  set of even $\xi$'s, $d(\xi)$ is even. Then a quadruple of type \lan even,
  odd, even, odd\ran\ has $d$-image of type \lan even, even, even, even\ran.
  Again this is impossible.
Similarly, there are no uncountable many even $\xi$'s with $d(\xi)$ odd.
 
  So there can be at most countably many changes of parity. In fact even if a
  single odd $\xi_1$ is with even $d(\xi_1)$, we get a contradiction by
  finding even $\xi_2,\xi_3,\xi_4$ so that $\langle \xi_1\ldots \xi_4\rangle$ is
  of type \lan odd, even, even, even\ran, but its $d$-image is of type \lan
  even, even, even, even\ran. Likewise, there is no even $\xi$ with odd
  $d(\xi)$.
 
  \itm{f} $d$ {\em is order preserving}. Since $d$ is one-to-one, it is order
  preserving on an uncountable set. First we prove that if $\xi_1<\xi_2$
  is of type \lan even, odd\ran\ or of type \lan odd, even\ran, then
  $d(\xi_1)<d(\xi_2)$. Assume this is not the case, and
  $d(\xi_1)>d(\xi_2)$. Then find $\xi_3,\xi_4$ such that
  $\seqn{\xi_1}{\xi_4}$ is of type $t_1=$\lan even, odd, even, odd\ran\
  or of type $t_2=$\lan odd, even, even, odd\ran\ and such that only $\xi_1$
  and $\xi_2$ change places.
 
  Then, since the first two coordinates change their place, the $d$-image of
  this quadruple is (respectively)
 of type $t_2$ or $t_1$. But since $t_1$ and $t_2$ are in
  different columns, a contradiction is derived.
 
  To see now that $d$ is order preserving on {\em any} pair, pick
  $\xi_1<\xi_2$ with the same parity. Then, since $\xi_1+1$ is of the
  other parity, $d(\xi_1)<d(\xi_1+1)<d(\xi_2)$.
 
  \itm{g}{\em $\gd$ is limit iff $d(\gd)$ is limit}. \\ 
This follows from
  the fact that \lan limit, even, odd, odd\ran\ is of type Suslin, while
  \lan even but not limit, even, odd, odd\ran\ is of type special.
 
  \itm{h} {\em For every $\xi$, $d(\xi+1)=d(\xi)+1$.} This is a
  consequence of the assumption that $\langle \xi,\xi+1,\mbox{odd,
  odd}\rangle$ is of type Suslin, while $\langle d(\xi),d(\xi+1),\mbox{ odd,
  odd}\rangle$ is of type special in case $d(\xi)+1<d(\xi +1)$.
 
  \itm{i}  $d$ {\em is onto} $\gw_1$. Since $d$ preserves the order, $I$ is
  well-ordered as well.
  Now, since for limit $\gd\in I$, $d(\gd)$ is limit, $d$ maps the block
  $[\gd,\gd+\gw)$ onto the block $[d(\gd),d(\gd)+\gw)$; thus the order-type
  of $\gd$ is the same as the order-type of $d(\gd)$ (use $E_\gd$). Hence
  $\gd=d(\gd)$ and $d$ is onto.
  \ed
 
  So we have concluded that $d$ is the identity, and thence the Unique Pattern
  Theorem.
  \section{ The encoding scheme}
  \label{s8}

  We now have all the ingredients for the main result---to show how to encode
  subsets of $\gw_1$ with simple patterns of Suslin sequences. The
  Magidor-Malitz quantifiers provide a concise way to describe our result.
 
  Recall that the quantified formula $Q x,y\ \gvp(x,y)$ holds in a structure
  $M$ if there is a set $\SA\subseteq M$ of cardinality $\aleph_1$, and for
  any two distinct $x,y\in \SA$, $M$ satisfies $\gvp(x,y)$. (Magidor and Malitz
  [1977]) 
 
  Let us see what can be stated in this language. We may say that the set
  $\SA$ (a unary predicate) is uncountable, simply by stating $Q x,y
  (\SA(x))$. To say that a linear order relation $\prec$ is $\gw_1$-like
   we just say that
  it is uncountable, and 
  any initial segment 
   $\setm{x}{x\prec y}$ is countable (and the obvious first-order
properties).
 
  Let us accept a slightly freer notion of $\gw_1$-trees: that of $\gw_1$-like
  trees. These are trees with set of levels, not $\gw_1$, but $\gw_1$-like.
The predecessors of a node in an $\go_1$-like tree form a countable chain
which is not necessarily well-ordered. Since an $\go_1$-like order embedds
$\go_1$, any $\go_1$-like tree contains an $\go_1$ tree.

  The notions of $\gw_1$-like Suslin trees, and $\gw_1$-like special trees can
  be defined and characterized in the Magidor-Malitz logic. There is a
  sentence $\sigma$ (in the language containing a binary relation $<$)
   such that $\ST\models \sigma$ iff $(\ST,<)$ is an
  $\gw_1$-like Suslin tree. $\sigma$ will simply state that $\ST$ is
 an uncountable tree (with obvious properties), and there is no
  uncountable set of pairwise incomparable  nodes in $\ST$. 

Going one more
  step, we describe now a sentence, $\gvp$, which holds true only in simple
  $\gw_1$-like sequences of $\gw_1$-like Suslin trees
  $\sseqn{\ST^\zeta}{\zeta\in I}$ where $I$ is an $\gw_1$-like ordering $\prec$.
  For this we may have to introduce a one-place predicate symbol $I$, and a
  two-place predicate $\ST(a,i)$ which, for a particular $i\in I$, describes
  the tree $\ST^i$. There is need also for
   a function which
  specializes those products which must be specialized. Giving more
  information on how $\gvp$ looks may annoy the reader who can find these
  details for herself; so we stop and state our theorem.
 
  \begin{theorem}[Encoding Theorem] 
\label{t8.1}
There is a sentence $\psi$ in the
  Magidor-Malitz logic which contains, besides the symbols $\prec$ etc.
  described above, a one place predicate $P(x)$ such that the following holds:
  Assuming $\diamond_{\gw_{1}}$, for any
  $X\subseteq \gw_1$, 
  \ben
  \item There is a model $M\models \psi$ for which $I^M=\gw_1$ and $P^N=X$.
  \item Assume $2^{\aleph_1}=\aleph_2$.
There is a generic extension of the universe which adds
  no new countable sets and collapses no cardinals, and such that in this
  extension:
If $N$ is any model satisfying $\psi$, then $I^N$ has order-type
  $\gw_1$, and (identifying $I^N$ with $\gw_1$)
  $$
  P^N=X
  $$
\een
  \end{theorem}
 
  \proof 
Let us first remark that the assumptions $\diamond_{\gw_1}$ and
 $2^{\aleph_1}=\aleph_2$ are not crucial, since these assumptions can be
obtained with a forcing of size $2^{2^{\aleph_1}}$.

The sentence $\psi$ is the simple-sequence sentence $\gvp$
partially described above with the addition
  of
  $$
  \A \zeta\in I(\langle 3,5,6,\zeta\rangle\in\st\mbox{ iff }P(\zeta))
  $$
 
  Given $X\subseteq\gw_1$ (assume that $X$ contains only ordinals $>6$), let
  $\langle \st,\spa\rangle$ be a simple pattern such that
  $X=\setm{\zeta\in\gw_1}{\langle 3,5,6,\zeta\rangle\in\st}$. Then use
  $\Diamond_{\gw_{1}}$ to construct an $\gw_1$-sequence
  $\CT=\sseqn{\ST^\zeta}{\zeta\in\gw_1}$ of Suslin trees which has the pattern
  $\langle\st,\spa\rangle$. This takes care of {\it 1}.

 The generic extension is a countable support
  iteration of specializing posets which keeps every derived tree of
  $\setm{\ST^d}{d\in\st}$ Suslin but specializes any \Aron\ tree, when possible.
  As we showed (Section ~\ref{s6}),
 an iteration of \length\ $\aleph_2$ suffices to make $\CT$ a
  primal sequence: Any \Aron\ tree $\SA$ in the extension is either special or
  it contains a club-embedding of a derived tree of some $\ST^d$, $d\in\st$.
   Now we have the required uniqueness.
  Let $N$ be a model of $\psi$. Then $I=I^N$ is an $\gw_1$-like ordering, and
  for each $i\in I$, a Suslin tree $\SA^i$ can be reconstituted from the set
  of $a$'s for which $\ST(a,i)$ holds in $N$. The sequence $\CA=
\sseqn{\SA^i}{i\in
  I}$ is simple and with pattern $\pair{\st^N}{\spa^N}$. The Unique Pattern
  Theorem ~\ref{t7.4}
can now be applied to $\CT$ and $\CA$ to yields $I^N\approx\gw_1$. And
  $\st=\st^N$. That is, the two sequences have the same pattern. From this it
  follows that $X=P^N$.
\subsection{The complete proof of Theorem A}
To obtain the $\Delta^2_2$ encoding of any subset of $\R$, Theorem
~\ref{t8.1} is sufficient. However, the statement of Theorem A in the
Preface is neater because it provides a categorical sentence, while Theorem
~\ref{t8.1} only establishes the uniqueness of $P^N$.
The sentence $\psi$ described above contains predicates and function
symbols other than $P$, and they too must be encoded by the unique pattern
sequence of trees.

To prove Theorem A, we ``catch our tail'' in the following way. 
Not only
the tree sequence $\langle \ST^\xi \rest \ga +1 \mid \xi \in \ga \rangle$
is constructed inductively, but so is the simple Suslin-special pattern
itself. More precisely, at the limit $\ga$ stage, we encode the countable
structure so far defined (the trees and $P$ etc.)
 as a subset of the ordinal-interval $(\ga, \ga+\gw
)$, and we put $\langle 3,5,6,\zeta \rangle \in \st$ for $\xi \in
(\ga, \ga+\gw)$ so that it encodes that countable structure. The
categorical sentence $\psi'$ tells us this fact as well.  The proof
continues just as before: any model $M$ for $\psi'$
determine a Suslin-special pattern which must be the unique such simple
pattern, but it determines in turn $M$, and hence the uniqueness of $M$.
\newpage

\section*{References}

\noindent

\rskip
{\sf U.Abraham and S. Shelah}
\ipar{1985}
{\em Isomorphism Types of Aronszajn Trees},
Israel J. of Math., Vol. 50, pp. 75-113.

\rskip
{\sf K. J. Devlin and H. Johnsbr\aa ten}
\ipar{1974} {The Souslin Problem},
Lecture Notes in Math. {\bf 405}, Springer-Verlag, New-York.

\rskip
{\sf T. Jech}
\ipar{1978}
{\em Set Theory}
Academic Press, New-York.

\rskip
{\sf R. B. Jensen and H. Johnsbr\aa ten}
\ipar{1974}
{\em A New Construction of a Non-Constructible $\Delta^1_3$ Subset of
$\gga$}, Fund. Math. 81, pp.279-290.

\rskip
{\sf M. Magidor and J. Malitz}
\ipar{1977}
{\em Compact Extensions of $L(Q)$},
Annals of Mathematical Logic, Vol. 11, pp. 217-261.

\rskip
{\sf S. Shelah}
\ipar{1982}
{\em Proper Forcing},
Lecture Notes in Math. {\bf 940}, Springer-Verlag, Berlin, Heidelberg,
New-York.
\ipar{1992}
{\em Proper and Improper Forcing}, (relevant preprints may be obtained from
the author).

\rskip
{\sf S. Shelah and H. Woodin}
\ipar{1990}
{\em Large cardinals imply every reasonably definable set is measurable},
Israel J. Math, 70, pp.381-394.

\rskip
{\sf S. Todor\v{c}evi\'{c}}
\ipar{1984}
{\em Trees and Linearly Ordered Sets},
in {\em Handbook of Set-Theoretic Topology} (K. Kunen and J. E. Vaughan,
eds), North-Holland, Amsterdam, pp. 235-293.
\rskip
{\sf H. Woodin}
\ipar{.}
{\em Large cardinals and determinacy}. in preparation.
\end{document}